\documentclass[11pt]{article} 

\newcommand{\doctype}{IEEE}

\usepackage{ifthen}

\long\def\comment#1{}

\newcommand{\myline}{1.5}

\newcommand{\citep}{\cite}

\usepackage{fullpage}

\usepackage{amsfonts,graphicx}
\usepackage{amsfonts,amsmath,amssymb,amsthm}
\usepackage{epic,eepic,eepicemu}
\usepackage{epsf}
\usepackage{graphics}
\usepackage{psfrag}
\usepackage{times}
\usepackage{epsfig}
\newtheorem{theorem}{Theorem}
\newtheorem{proposition}[theorem]{Proposition}
\newtheorem{lemma}[theorem]{Lemma}
\newtheorem{definition}[theorem]{Definition}

\newtheorem{cor}[theorem]{Corollary}

\makeatletter
\def\@cite#1#2{[\if@tempswa #2 \fi #1]}
\makeatother

\def\E{{\mathbb E}}

\def\indicator{{\mathbb I}}

\def\var{{\rm Var}}

\newcommand{\conv}{\ensuremath{\operatorname{conv}}}

\newcommand{\dims}{\ensuremath{d}}
\newcommand{\ind}{\ensuremath{k}}

\newcommand{\Uspace}{\ensuremath{\mathcal{U}}}
\newcommand{\Umax}{\ensuremath{K}}
\newcommand{\firstPrior}{\ensuremath{\pi^1}}
\newcommand{\secondPrior}{\ensuremath{\pi^0}}
\newcommand{\timestep}{\ensuremath{n}}
\newcommand{\minusOne}{\ensuremath{0}}

\newcommand{\asymKL}[1]{\ensuremath{\mykull^0_{#1}}}
\newcommand{\asymLK}[1]{\ensuremath{\mykull^1_{#1}}}
\newcommand{\extKL}[2]{\ensuremath{\tilde{D}^0(#1,#2)}}
\newcommand{\extLK}[2]{\ensuremath{\tilde{D}^1(#1,#2)}}
\newcommand{\extKLplain}{\ensuremath{\tilde{D}^0}}
\newcommand{\extLKplain}{\ensuremath{\tilde{D}^1}}

\newcommand{\derivone}[2]{\ensuremath{#1_{#2}}}
\newcommand{\derivtwo}[3]{\ensuremath{#1_{#2 #3}}}

\newcommand{\defn}{\ensuremath{:=}}

\newcommand{\real}{\ensuremath{\mathbb{R}}}

\newcommand{\Xspace}{\ensuremath{\mathcal{X}}}

\newcommand{\widgraph}[2]{\includegraphics[keepaspectratio,width=#1]{#2}}

\newcommand{\vtiny}{\ensuremath{\vspace*{.01in}}}

\newcommand{\myparagraph}[1]{{\bf{#1}}}

\newcommand{\mysection}[1]{\vtiny \\ {\bf{#1}}}

\makeatletter
\long\def\@makecaption#1#2{
        \vskip 0.8ex
        \setbox\@tempboxa\hbox{\small {\bf #1:} #2}
        \parindent 1.5em  
        \dimen0=\hsize
        \advance\dimen0 by -3em
        \ifdim \wd\@tempboxa >\dimen0
                \hbox to \hsize{
                        \parindent 0em
                        \hfil 
                        \parbox{\dimen0}{\def\baselinestretch{0.96}\small
                                {\bf #1.} #2
                                } 
                        \hfil}
        \else \hbox to \hsize{\hfil \box\@tempboxa \hfil}
        \fi
        }
\makeatother

\makeatletter
\long\def\barenote#1{
    \insert\footins{\footnotesize
    \interlinepenalty\interfootnotelinepenalty 
    \splittopskip\footnotesep
    \splitmaxdepth \dp\strutbox \floatingpenalty \@MM
    \hsize\columnwidth \@parboxrestore
    {\rule{\z@}{\footnotesep}\ignorespaces
      #1\strut}}}
\makeatother

\newcommand{\Prob}{\ensuremath{\mathbb{P}}}
\newcommand{\Cprob}[2]{\ensuremath{\Prob(\, #1 \, | \, #2 \,)}}

\newcommand{\mybeginproof}{\noindent \emph{Proof: $\;$}}
\newcommand{\myendproof}{\hfill \qed}

\newcommand{\fq}{f_{\phi_\timestep}} \newcommand{\fqs}{f_{\phi}}

\newcommand{\up}{\ensuremath{M}}
\newcommand{\mykullqs}{\ensuremath{\mykull_{\phi}}}

\newcommand{\Jfunc}{\ensuremath{J}}
\newcommand{\Gfunc}{\ensuremath{G}}
\newcommand{\Jtil}{\ensuremath{\widetilde{J}}}

\newcommand{\Jfuncphi}{\ensuremath{\Jfunc_\phi}}
\newcommand{\Gfuncphi}{\ensuremath{\Gfunc_\phi}}

\newcommand{\ess}{\ensuremath{\operatorname{ess}}}

\newcommand{\mykull}{\ensuremath{D}}

\begin{document}

\comment{
\ifthenelse{\equal{\doctype}{IEEE}}
{
\typeout{---------> No formatting modification for IEEE}
\typeout{}
}
{
\typeout{------> Setting Stat. tech report formatting}
\typeout{}
\oddsidemargin .00in
\evensidemargin .00in
\marginparwidth 0.07 true in
\addtolength{\headsep}{0.25in}
\textheight 8.53 true in
\textwidth 6.53 true in
\clubpenalty=10000
}
}

\title{On optimal quantization rules for some problems in sequential
decentralized detection}

\author{XuanLong Nguyen, Martin J. Wainwright and Michael I. Jordan}

\comment{
\ifthenelse{\equal{\doctype}{IEEE}}
{
\author{ XuanLong Nguyen$^{1}$, Martin J. Wainwright$^{2}$ and 
Michael I. Jordan$^{2}$ \\ 
$^1$ SAMSI and Dept. of Statistical Science, Duke University\\
$^2$ Dept. of EECS and Dept. of Statistics, UC Berkeley \\
\texttt{\small xuanlong.nguyen@gmail.com;\{wainwrig,jordan\}@eecs.berkeley.edu} \\
}}
{
\author{ XuanLong Nguyen\\ Department of EECS \\ 
University of California, Berkeley\\ 
\texttt{\small xuanlong@cs.berkeley.edu} \\[.2cm]
\and
Martin J. Wainwright \\ Department of Statistics and Department of EECS \\
University of California, Berkeley\\ \texttt{\small
wainwrig@stat.berkeley.edu} \\[.2cm]
\and
Michael I. Jordan \\ Department of Statistics and Department of EECS\\ 
University of California, Berkeley\\ \texttt{\small
jordan@stat.berkeley.edu}\\ 
}}
\vspace{.2cm}
}

\date{}

\maketitle


\comment{
\ifthenelse{\equal{\doctype}{IEEE}}
{\typeout{}}
{ \begin{center}
Technical Report 708
Department of Statistics \\
University of California, Berkeley\\
\bigskip
\end{center}
}

\ifthenelse{\equal{\doctype}{IEEE}}
{\renewcommand{\baselinestretch}{\myline}}
{ }
}

\begin{abstract}
We consider the design of systems for sequential decentralized
detection, a problem that entails several interdependent choices: the
choice of a stopping rule (specifying the sample size), a global
decision function (a choice between two competing hypotheses), and a
set of quantization rules (the local decisions on the basis of which
the global decision is made).  This paper addresses an open problem of
whether in the Bayesian formulation of sequential decentralized
detection, optimal local decision functions can be found within the
class of stationary rules.  We develop an asymptotic approximation to
the optimal cost of stationary quantization rules and exploit this
approximation to show that stationary quantizers are not optimal in
a broad class of settings.  We also consider the class of blockwise 
stationary quantizers, and show that asymptotically optimal quantizers are
likelihood-based threshold rules.\footnote{
XuanLong Nguyen is with the Statistical and Applied
Mathematical Sciences Institute and Duke University's
Department of Statistical Science. \\
Martin J. Wainwright is with the Department of Statistics
and Department of Electrical Engineering and 
Computer Sciences, University of California, Berkeley. \\
Michael I. Jordan is with the Department of Statistics
and Department of Electrical Engineering and 
Computer Sciences, University of California, Berkeley. \\
This work was presented in
part at the International Symposium on Information Theory, July 2006,
Seattle, WA.}

\end{abstract}

\noindent {\bf{Keywords:}} decentralized
detection; decision-making under constraints; experimental design; 
hypothesis testing; quantizer design; sequential detection.

\section{Introduction}

Detection is a classical discrimination or hypothesis-testing problem,
in which observations $\{X_1, X_2, \ldots\}$ are assumed to be drawn
i.i.d. from the (multivariate) conditional distribution
\mbox{$\Cprob{\cdot}{H}$} and the goal is to infer the value of the
random variable $H$, which takes values in $\{0,1\}$.  In a typical
engineering application, the case $\{H=1\}$ represents the presence of
some target to be detected, whereas $\{H=0\}$ represents its absence.
Placing this problem in a communication-theoretic context, a
\emph{decentralized detection} problem is a hypothesis-testing problem
in which the decision-maker is not given access to the raw data points
$X_n$, but instead must infer $H$ based only on the output of a set of
\emph{quantization rules} or \emph{local decision functions}, say
$\{U_n = \phi_n(X_n)\}$, which map the raw data to quantized values.
This basic problem of decentralized detection has been studied
extensively for several decades~\cite{Tenney81,Tsitsiklis86,Chen05};
see the overview
papers~\cite{Tsitsiklis93,Viswanathan97,Blum97,Chen06} and references
therein for more background.  Of interest in this paper is the
extension to an-online setting: more specifically, the
\emph{sequential decentralized detection}
problem~\cite{Tsitsiklis86,Veeravalli99,Mei-thesis} involves a data
sequence, $\{X_1, X_2, \ldots\}$, and a corresponding sequence of
summary statistics, $\{U_1, U_2, \ldots\}$, determined by a sequence
of local decision rules $\{\phi_1, \phi_2, \ldots \}$.  The goal is to
design both the local decision functions and to specify a global
decision rule so as to predict $H$ in a manner that optimally trades
off accuracy and delay.  In short, the sequential decentralized
detection problem is the communication-constrained extension of
classical formulation of sequential centralized decision-making
problems (see, e.g.,~\cite{Chernoff72,Shiryayev78,Lai01}) to the
decentralized setting.

\comment{ In particular, in Case E, the local sensors are provided
with memory and with feedback from the global decision-maker (also
known as the \emph{fusion center}), so that each sensor has available
to it the current data, $X_n$, as well as all of the summary
statistics from all of the other local sensors.  In other words, each
local sensor has the same snapshot of past state as the fusion center;
this is an instance of a so-called ``quasi-classical information
structure''~\cite{Ho80} for which dynamic programming (DP)
characterizations of the optimal decision functions are available.
Veeravalli et al.~\cite{Veeravalli93} exploit this fact to show that
the decentralized case has much in common with the centralized case,
in particular that likelihood ratio tests are optimal local decision
functions at the sensors and that a variant of a sequential
probability ratio test is optimal for the decision-maker.
}

In setting up a general framework for studying sequential
decentralized problems, Veeravalli et al.~\cite{Veeravalli93} defined
five problems, denoted ``Case A'' through ``Case E,'' distinguished
from one another by the amount of information available to the local
sensors.  In applications such as power-constrained sensor networks,
one cannot assume that the decision-maker and sensors can communicate
over a high-bandwidth channel, nor that the sensors have unbounded
memory.  Most suited to this perspective---and the focus of this
paper---is Case A, in which the local decisions are of the simplified
form $\phi_n(X_n)$; i.e., neither local memory nor feedback are
assumed to be available.  Noting that Case A is not amenable to
dynamic programming and hence presumably intractable, Veeravalli et
al.~\cite{Veeravalli93} suggested restricting the analysis to the
class of \emph{stationary} local decision functions; i.e., local
decision functions $\phi_n$ that are independent of $n$.  They
conjectured that stationary decision functions might actually be
optimal in the setting of Case A (given the intuitive symmetry and
high degree of independence of the problem in this case), even though
it is not possible to verify this optimality via DP arguments.  This
conjecture has remained open since it was first posed by Veeravalli et
al.~\cite{Veeravalli93,Veeravalli99}.

The main contribution of this paper is to resolve this question by
showing that stationary decision functions are, in fact, \emph{not}
optimal for decentralized problems of type A.  Our argument is based
on an asymptotic characterization of the optimal Bayesian risk as the
cost per sample goes to zero.  In this asymptotic regime, the optimal
cost can be expressed as a simple function of priors and
Kullback-Leibler (KL) divergences.  This characterization allows us to
construct counterexamples to the stationarity conjecture, both in an
exact and an asymptotic setting.  In the latter setting, we present a
class of problems in which there always exists a range of prior
probabilities for which stationary strategies, either deterministic or
randomized, are suboptimal.  We note in passing that an 
intuition for the source of the suboptimality is easily 
provided---it is due to the asymmetry of the KL divergence.

It is well known that optimal quantizers when unrestricted are
necessarily likelihood-based threshold rules~\cite{Tsitsiklis86}.  Our
counterexamples and analysis imply that optimal thresholds are not
generally stationary (i.e., the threshold may differ from sample to
sample).  We also provide a partial converse to this result:
specifically, if we restrict ourselves to stationary (or blockwise
stationary) quantizer designs, then there exists an optimal design
that is a deterministic threshold rule based on the likelihood 
ratio. We prove this result by establishing a quasiconcavity 
result for the asymptotically optimal cost function.  

It is worth highlighting several limitations in our results.
For the suboptimality of stationary quantizers, our analysis is 
applicable only to finite classes of deterministic quantizers and 
their convex hull of randomized quantizers, and under the 
assumption that the likelihood ratio of the two hypotheses are 
bounded from both above and below.  Such assumptions certainly hold 
for arbitrary discrete distributions with finite support. 
It remains an open problem to consider more general classes
of distributions.  For the likelihood-ratio characterization result, 
our proof works only for the (possibly infinite) classes of 
deterministic quantizers with arbitrary output alphabets, as well 
as for the class of randomized quantizers with binary outputs.  
We conjecture that the same result holds more generally for 
randomized quantizers with arbitrary output alphabets.

The remainder of this paper is organized as follows.  We begin in
Section~\ref{SecBackground} with background on the Bayesian
formulation of sequential detection problems, and Wald's
approximation.  Section~\ref{SecCharOpt} provides a simple asymptotic
approximation of the optimal cost that underlies our main analysis in
Section~\ref{SecSuboptStat}.  In Section~\ref{SecAsympBlock}, we
establish the existence of optimal decision rules that are
likelihood-based threshold rules, under the restriction to blockwise
stationarity.  We conclude with a discussion in
Section~\ref{SecDiscussion}.


\section{Background}
\label{SecBackground}

This section provides background on the Bayesian formulation of
sequential (centralized) detection problems.  Of particular use in our
subsequent analysis is Wald's approximation of the cost of optimal
sequential test.

Let $\Prob_0$ and $\Prob_1$ represent the distribution of $X$, when
conditioned on $\{H=0\}$ and $\{H=1\}$ respectively.  Assume that
$\Prob_0$ and $\Prob_1$ are absolutely continuous with respect to one
another.  We use $f^0(x)$ and $f^1(x)$ to denote the respective
density functions with respect to some dominating measure (e.g.,
Lebesgue for continuous variables, or counting measure for
discrete-valued variables).

Our focus is the Bayesian formulation of the sequential detection
problem~\cite{Shiryayev78,Veeravalli99}; accordingly, we let $\pi^1 =
\Prob(H=1)$ and $\pi^0 = \Prob(H=0)$ denote the prior probabilities of
the two hypotheses. Let $X_1, X_2,\ldots$ be a sequence of conditionally
i.i.d. realizations of $X$. A sequential decision rule consists of a
\emph{stopping time} $N$ defined with respect to the sigma field
$\sigma(X_1,\ldots,X_N)$, and a decision function $\gamma$ measurable
with respect to $\sigma(X_1,\ldots,X_N)$.  The cost function is the
expectation of a weighted sum of the sample size $N$ and the probability 
of incorrect decision---namely
\begin{eqnarray}
\label{EqnDefnLoss}
J(N,\gamma) \defn \E \big \{ cN + \indicator[\gamma(X_1,\ldots,X_N)
\neq H] \big \},
\end{eqnarray}
where $c > 0$ is the incremental cost of each sample.  The overall
goal is to choose the pair $(N, \gamma)$ so as to minimize the
expected loss~\eqref{EqnDefnLoss}.

It is well known that the optimal solution of the sequential decision
problem can be characterized recursively using dynamic programming
(DP)
arguments~\cite{Arrowetal49,Waldetal48,Shiryayev78,Bertsekas_dyn1}.
Although useful in classical (centralized) sequential detection, the
DP approach is not always straightforward to apply to
\emph{decentralized} versions of sequential
detection~\cite{Veeravalli99}.  In the remainder of this section, we
describe an asymptotic approximation of the optimal sequential cost,
originally due to Wald (cf.~\cite{Siegmund85}), valid as $c
\rightarrow 0$.  To sketch out Wald's approximation, we begin by
noting the optimal stopping rule for the cost
function~\eqref{EqnDefnLoss} takes the form
\begin{equation}
\label{EqnDefnStopping}
N = \inf \big \{ n \geq 1 \; \big | \; L_n(X_1, \ldots, X_n) \defn
\sum_{i=1}^n \log \frac{f^1(X_i)}{f^0(X_i)} \notin (a,b) \big \},
\end{equation}
for some real numbers $a < b$.  Given this stopping rule, the optimal
decision function has the form
\begin{equation}
\label{EqnDefnFinalDec}
\gamma(L_N) = \begin{cases} 1 & \mbox{if $L_N \geq b$}, \\
0 & \mbox{if $L_N \leq a$.}
			   \end{cases}
\end{equation}

Consider the two types of error:
\begin{eqnarray*}
\alpha & = & \Prob_0(\gamma(L_N) \neq H ) = \Prob_0(L_N \geq b) \\
\beta & = & \Prob_1(\gamma(L_N) \neq H) = \Prob_1(L_N \leq a).
\end{eqnarray*}
As $c\rightarrow 0$, it can be shown that the optimal choice of $a$
and $b$ satisfies $a\rightarrow -\infty, b\rightarrow \infty$, and the
corresponding $\alpha,\beta$ satisfy $\alpha + \beta \rightarrow 0$.
Ignoring the overshoot of $L_N$ upon the optimal stopping time $N$
(i.e., instead assuming $L_N$ attains precisely the value $a$ or $b$)
we can express $a$, $b$, $\E N$ and the cost function $J$ in terms of
$\alpha$ and $\beta$ as follows~\cite{Wald47}:
\begin{eqnarray}
\label{Approx-ab}
a \approx a(\alpha, \beta) \defn \log \frac{\beta}{1-\alpha} \;
& \mbox{and} & \; b \approx b(\alpha, \beta) \defn \log
\frac{1-\beta}{\alpha} \\
\label{Approx-KL}
\E_0[L_N]  \approx  (1-\alpha)a + \alpha b \;
& \mbox{and} & \;
\E_1[L_N] \approx  (1-\beta) b + \beta a 
\end{eqnarray}

Now define the Kullback-Leibler divergences 
\begin{equation}
\mykull^1 = \E_1[\log \frac{f^1(X_1)}{f^0(X_1)}] = D(f^1||f^0), \qquad
\mbox{and} \qquad \mykull^0 = -\E_0[\log \frac{f^1(X_1)}{f^0(X_1)}] =
D(f^0||f^1).
\end{equation}
With a slight abuse of notation, we shall also use $D(\alpha,\beta)$
to denote a function in $[0,1]^2 \rightarrow \real$ such that:
\[D(\alpha,\beta) \defn \alpha \log\frac{\alpha}{\beta} + (1-\alpha)\log
\frac{1-\alpha}{1-\beta}.\] With the above approximations, the cost
function $\Jfunc$ of the decision rule based on envelopes $a$ and $b$
can be written as
\begin{eqnarray}
\label{EqnDefnJcost}
%
\Jfunc & = & \pi^1 \E_1 (cN + \indicator[L_N \leq a]) +
\pi^0 \E_0 (cN + \indicator[L_N \geq b]) \nonumber \\
& = & c \pi^1 \frac{\E_1 L_N}{\mykull^1} + c \pi^0 \frac{\E_0 L_N}{-\mykull^0}
 + \pi^0 \alpha + \pi^1 \beta, \qquad \\
\label{EqnDefnGfunc}
& \approx & c \secondPrior
\frac{D(\alpha,1-\beta)}{\mykull^0} + c \firstPrior
\frac{D(1-\beta,\alpha)}{\mykull^1} 
+ \secondPrior \alpha + \firstPrior \beta,
\end{eqnarray}
where the third line follows from Wald's equation~\cite{Wald47}. 
Let $\Jtil(\alpha,\beta)$ denote the approximation~\eqref{EqnDefnGfunc} of $J$.

Let $\Jfunc^*$ denote the cost of an optimal sequential test, i.e.,
\begin{equation}
\label{EqnDefnOptJ}
\Jfunc^* = \inf_{a,b} J.
\end{equation}
\comment{
Since $\alpha + \beta \rightarrow 0$, $D(1-\beta,\alpha) = \log (1/\alpha) + o(1)$, and
$D(1-\alpha,\beta) = \log(1/\beta) + o(1)$.  We approximate
$J^*$ by minimizing $\Jfunc$ over $\alpha$ and $\beta$. 
The minimum is achieved at
$\alpha^* = \frac{c\pi^1}{\mykullqs^1\pi^0}$ and $\beta^* =
\frac{c\pi^0}{\mykullqs^0\pi^1}$, yielding:
\begin{eqnarray}
\Jfunc^* & \approx & \inf_{\alpha,\beta} \biggr \{ \pi^0 \alpha + \pi^1 \beta +
c\pi^0\frac{\log(1/\beta)}{\mykull^0} +
c\pi^1\frac{\log(1/\alpha)}{\mykull^1} \biggr \} + o(c) \nonumber \\
\label{Eqn-Approx}
&\approx & (\frac{\pi^0}{\mykull^0} + \frac{\pi^1}{\mykull^1}) c\log c^{-1} + O(c).
\end{eqnarray}
} A useful result due Chernoff~\cite{Chernoff59} states that under
certain assumption (to be elaborated in the next section), $\Jfunc^*$
has the following form:
\begin{eqnarray}
\Jfunc^* 
\label{Eqn-Approx}
&\approx & (\frac{\pi^0}{\mykull^0} + \frac{\pi^1}{\mykull^1}) c\log c^{-1} (1 + o(1)).
\end{eqnarray}

\section{Characterization of optimal stationary quantizers}
\label{SecCharOpt}

Turning now to the decentralized setting, the primary challenge lies
in the design of the quantization rules $\phi_\timestep$ applied to
data $X_n$.  When $X_n$ is univariate, a deterministic quantization
rule $\phi_\timestep$ is a function that maps $\Xspace$ to the
discrete space $\Uspace = \{0,\ldots, \Umax - 1\}$ for some natural
number $\Umax$. For multivariate $X_n$ with $\dims$ dimensions arising
in the multiple sensor setting, a deterministic quantizer $\phi_n$
is defined as a mapping from the $\dims$-dimensional product space $\Xspace$
to $\Uspace = \{0, \ldots, \Umax - 1\}^\dims$. In the decentralized
problem defined as Case A by Veeravalli et al.~\cite{Veeravalli93},
the function $\phi_n$ is composed of $\dims$ separate quantizer
functions, one each for each dimension.  A randomized quantizer $\phi_n$
is obtained by placing a distribution over the space of deterministic
quantizers.

Any fixed set of quantization rules $\phi_\timestep$ yields a sequence
of compressed data $U_\timestep=\phi_\timestep(X_\timestep)$, to which
the classical theory can be applied.  We are thus interested in
choosing quantization rules $\phi_1,\phi_2,\ldots$ so that the error
resulting from applying the optimal sequential test to the sequence 
of statistics $U_1,U_2,\ldots$ is minimized over some space
$\Phi$ of quantization rules. For a given quantizer $\phi_n$ we use
\begin{eqnarray*}
\fq^i(u) & \defn & \Prob_i(\phi_\timestep(X_\timestep) = u), \qquad
\mbox{for} \quad i = 0,1,
%
\end{eqnarray*}
to denote the distributions of the compressed data, conditioned on the
hypothesis. In general, when randomized quantizers are allowed, the
vector $(\fq^0(.),\fq^1(.))$ ranges over a convex set, denoted $\conv
\Phi$, whose extreme points correspond to deterministic quantizers
based on likelihood ratio threshold
rules~\cite{Tsitsiklis93-extremal}.

We say that a quantizer design is \emph{stationary} if the rule
$\phi_n$ is independent of $n$; in this case, we simplify the notation
to $\fqs^1$ and $\fqs^0$.  In addition, we define the KL divergences
$\mykullqs^1 \defn D(\fqs^1||\fqs^0)$ and $\mykullqs^0 \defn
D(\fqs^0||\fqs^1)$.  Moreover, let $\Jfuncphi$ and $\Jfuncphi^*$
denote the analogues of the functions $\Jfunc$ in
Eq.~\eqref{EqnDefnJcost} and $\Jfunc^*$ in Eq.~\eqref{EqnDefnOptJ}, 
respectively, 
defined using $\mykullqs^i$, for $i=0,1$.  In this scenario, the sequence 
of compressed data $U_1,\ldots, U_n,\ldots$ are drawn i.i.d. from either 
$\fqs^0$ or $\fqs^1$.  Thus we can use the approximation~\eqref{Eqn-Approx} 
to characterize the asymptotically optimal stationary quantizer design. 
This is stated formally in the lemma to follow.

We begin by stating the assumptions underlying the lemma.  For a given
class of quantizers $\Phi$, we assume that the Kullback-Leibler divergences
are uniformly bounded away from zero
\begin{eqnarray}
\label{Eqn-Ass1}
&& D(\fqs^1||\fqs^0) > 0, D(\fqs^0||\fqs^1) > 0 \;\mbox{for all}\;
\phi \in \Phi
\end{eqnarray}
and moreover that the variance of the log likelihood ratios are bounded
\begin{eqnarray}
\label{Eqn-Ass2}
&& \sup_{\phi\in \Phi} \var_{\fqs^1}\log(\fqs^1/\fqs^0) < \infty,
\qquad \mbox{and} \qquad \sup_{\phi \in \Phi}
\var_{\fqs^0}\log(\fqs^1/\fqs^0) < \infty.
\end{eqnarray}
\begin{lemma}
\label{PropApprox}
(a) Under assumptions~\eqref{Eqn-Ass1} and~\eqref{Eqn-Ass2}, the optimal
stationary cost takes the form
\begin{equation}
\label{EqnOptJfuncphi} 
\Jfuncphi^* = \biggr (\frac{\pi^0}{\mykullqs^0} + \frac{\pi^1}{\mykullqs^1}
\biggr ) \; c\log c^{-1} \; (1+ r_\phi)
\end{equation}
where $|r_\phi| = o(1)$ as $c\rightarrow 0$.

(b) If $\sup_{\phi \in \Phi} \max \{\log(\fqs^1/\fqs^0), \log(\fqs^0/\fqs^1)\} 
< M$ for some constant $M$, then~\eqref{EqnOptJfuncphi} holds 
with \\
$\sup_{\phi \in \Phi} |r_\phi| = o(1)$ as $c\rightarrow 0$.
\end{lemma}
\mybeginproof 
(a) This part is immediate from a combination of Theorems
1 and 2 of~Chernoff~\cite{Chernoff59}.
\comment{
We prove the lemma using results originally due to
Chernoff~\cite{Chernoff59}, restricted to a simple binary hypothesis
test between $\fqs^0$ and $\fqs^1$. By Theorem 1 from
Chernoff~\cite{Chernoff59}, under conditions~\eqref{Eqn-Ass1}
and~\eqref{Eqn-Ass2}, there is a sequential test $(N,\gamma)$ for
which:
\begin{eqnarray*}
\Jfunc^* \leq \Jfunc(N,\gamma) & = & \pi^0(\alpha + c\E_0 N) +
\pi^1(\beta + c\E_1 N) \\ & \leq & \pi^0 (1 + o(1))c\log
c^{-1}/\mykullqs^0 + \pi^1 (1 + o(1)) c\log c^{-1}/\mykullqs^1.
\end{eqnarray*}
But then the optimal test with the cost $\Jfunc^*$ (i.e., the
likelihood ratio based test) must satisfies that $\alpha + c\E_0 N =
O(c\log c^{-1})$ and $\beta + c\E_1 N = O(c\log c^{-1})$.  Theorem 2
of Chernoff~\cite{Chernoff59} implies that
\[\Jfunc^* \geq \biggr (\frac{\pi^0}{\mykullqs^0} + \frac{\pi^1}{\mykullqs^1}
\biggr ) (1+o(1))c\log c^{-1},\] concluding the proof.
}

(b) We begin by bounding the error in the
approximation~\eqref{EqnDefnGfunc}.  
By definition of the stopping time
$N$, we have either (i) $b \leq L_N \leq b + \up $ or (ii) $a - \up
\leq L_N \leq a$.  
\comment{Consider all realizations $u_1,\ldots,u_n$ for
which condition (i) holds; for any such sequence, we have
\begin{eqnarray*}
e^b P_0(u_1,\ldots,u_n) & \leq & P_1(u_1,\ldots,u_n) \\
e^{(b+\up)} P_0(u_1,\ldots,u_n) & \geq & P_1(u_1,\ldots,u_n).
\end{eqnarray*}
Taking a sum over all such realizations, using the definition of
$\alpha$ and $\beta$, and performing some algebra yields the
inequality}
By standard arguments due to Wald~\cite{Wald47}, it is simple
to obtain $e^b \alpha \, \leq \, 1-\beta \, \leq \, e^{b+\up}\alpha$,
or equivalently $b \leq b(\alpha,\beta) = \log \frac{1-\beta}
{\alpha} \leq b + \up$.
Similar reasoning for case (ii) yields $a - \up \leq a(\alpha,\beta)=\log
\frac{\beta}{1-\alpha} \leq a$. Now, note that
\[\E_0 L_N = \alpha \E_0[L_N | L_N \geq b] + (1-\alpha) \E_0[L_N|L_N \leq a].\]
Conditioning on the event $L_N \in [b, b+\up]$, we have $|L_N -
b(\alpha,\beta)| \leq M$. Similarly, conditioning on the event $L_N
\in [a-\up,a]$, we have $|L_N - b(\alpha,\beta)| \leq M$.  This yields
$|\E_0 L_N - (-D(\alpha,1-\beta))| \leq M$.  Similar reasoning yields
$|\E_1 L_N - D(1-\beta,\alpha)| \leq M$. Let $\Jtil_\phi(a,b)$
denote the approximation~\eqref{EqnDefnGfunc} of $J_\phi$.
We obtain:
\[|J_\phi - \Jtil_\phi(\alpha,\beta)| \leq 2cM.\]

Note that the approximation error bound is independent of $\phi$.
Thus, it suffices to establish the asymptotic
behavior~\eqref{EqnOptJfuncphi} for the quantity 
$\inf_{\alpha, \beta} \Jtil_\phi(\alpha, \beta)$, where the infimum is
taken over pairs of realizable error probabilities
$(\alpha,\beta)$. Moreover, we only need to consider the asymptotic
regime $\alpha+\beta \rightarrow 0$, since the error probabilities
$\alpha$ and $\beta$ vanish as $c \rightarrow 0$.  It is simple to see
that $D(1-\beta,\alpha) = \log (1/\alpha)(1 + o(1))$, and
$D(1-\alpha,\beta) = \log(1/\beta)(1 + o(1))$.  Hence, $\inf_{\alpha,
\beta} \Jtil_\phi(\alpha, \beta)$ can be expressed as
\begin{equation}
\label{EqnInf}
\inf_{\alpha,\beta} \biggr \{ \pi^0 \alpha + \pi^1 \beta +
c\pi^0\frac{\log(1/\beta)}{\mykullqs^0} +
c\pi^1\frac{\log(1/\alpha)}{\mykullqs^1} \biggr \}(1 + o(1)).
\end{equation}
This infimum, taken over all positive $(\alpha,\beta)$, is achieved at
\mbox{$\alpha^* = \frac{c\pi^1}{\mykullqs^1\pi^0}$} and \mbox{$\beta^* =
\frac{c\pi^0}{\mykullqs^0\pi^1}$}. 
Plugging the quantities $\alpha^*$
and $\beta^*$ into Eq.~\eqref{EqnInf} yields~\eqref{EqnOptJfuncphi}.
Note that the asymptotic quantity $o(1)$ in~\eqref{EqnOptJfuncphi}
is absolutely bounded by
$\alpha^* + \beta^* \rightarrow 0$ uniformly for all quantizer $\phi$,
because $\mykullqs^1$ and $\mykullqs^0$ are uniformly bounded away from 
zero due to the Lemma's assumption.

It remains to show that error probabilities $(\alpha^*,\beta^*)$
can be approximately realized by using a
sufficiently large threshold $b > 0$ and small threshold $a < 0$ while
incurring an approximation cost of order $O(c)$ uniformly for
all $\phi$.  Indeed, let us choose
thresholds $a'$ and $b'$ such that $e^{-(b'+M)}/2 \leq \alpha^* \leq
e^{-b'}$, and $e^{a'-M}/2 \leq \beta^* \leq e^{a'}$. Let $\alpha'$ and
$\beta'$ be the corresponding errors associated with these two
thresholds. As before, we also have $\alpha' \in
(e^{-(b'+M)}/2,e^{-b'})$ and $\beta' \in (e^{a'-M}/2,
e^{a'})$. Clearly, $|\alpha^* - \alpha'| \leq e^{-b'}(1-e^{-M}/2) =
O(\alpha^*) = O(c)$. Similarly, $|\beta^* - \beta'| = O(c)$.  By the
mean value theorem,
\[|\log(1/\alpha^*) - \log(1/\alpha')| \leq 
|\alpha^*-\alpha'|e^{b'+M} \leq 2e^{M}(1-e^{-M}/2) = O(1).\]
Similarly, $\log(1/\beta^*) - \log(1/\beta') = O(1)$.  Hence, the
approximation of $(\alpha^*,\beta^*)$ by the realizable
$(\alpha',\beta')$ incurs a cost at most $O(c)$. Furthermore,
the constant in the asymptotic bound $O(c)$ is independent
of quantizer $\phi \in \Phi$.

\comment{
Plugging the quantities $\alpha^*,\beta^*$
into equation~\eqref{EqnInf} yields 
\begin{subequations}
\begin{eqnarray}
\label{Eqn-approx}
\inf_{\alpha, \beta} \Gfuncphi(\alpha, \beta) & = & \biggr
(\frac{\pi^0}{\mykullqs^0} + \frac{\pi^1}{\mykullqs^1} \biggr )
c\log\frac{1}{c} + \biggr (\frac{\pi^0}{\mykullqs^0} +
\frac{\pi^1}{\mykullqs^1} + \frac{\pi^0}{\mykullqs^0}\log
\frac{\mykullqs^0\pi^1}{\pi^0} + \frac{\pi^1}{\mykullqs^1}\log
\frac{\mykullqs^1\pi^0}{\pi^1} \biggr ) c + O(c) \\ & = & \biggr
(\frac{\pi^0}{\mykullqs^0} + \frac{\pi^1}{\mykullqs^1} \biggr )
c\log\frac{1}{c} + O(c)
\end{eqnarray}
\end{subequations}
as claimed.
}

%

\myendproof

\noindent \myparagraph{Remarks:} 
\begin{enumerate}
\item If $\Phi$ is a finite class of quantizers, or
a convex hull of a finite class of quantizers,
the assumption in part b of Lemma~\ref{PropApprox} holds. 
It also holds in the case of discrete distributions and
continuous distributions with bounded support.
However, it would be interesting to relax this assumption so 
as to cover distributions with unbounded support.  

\item The preceding approximation of the optimal cost essentially
ignores the overshoot of the likelihood ratio $L_N$. While it is
possible to analyze this overshoot to obtain a finer approximation
(cf.~\cite{Lorden70,Siegmund85,Lai01,Poor94}), we see that this is not
needed for our purpose. Lemma~\ref{PropApprox} shows that given a
fixed prior $(\pi^0,\pi^1)$, among all stationary quantizer designs in
$\Phi$, $\phi$ is optimal for sufficiently small $c$ if \emph{and}
only if $\phi$ minimizes what we shall call the \emph{sequential cost
coefficient}:
\[G_\phi \defn \frac{\pi^0}{\mykullqs^0} + \frac{\pi^1}{\mykullqs^1}.\]

\item As a consequence of Lemma~\ref{Lem-Quasiconcave} to be proved in
the sequel, if we consider the class $\Phi$ of all binary randomized
quantizers, then sequential cost coefficient $G_\phi$ is a
quasiconcave function with respect to $(\fqs^0(.),\fqs^1(.))$.  (A
function $F$ is quasiconcave if and only if for any $\eta$, the level
set $\{F(x) \geq \eta \}$ is a convex set; see Boyd and
Vandenberghe~\cite{Boyd-book} for further background).  The minimum of
a quasiconcave function lies in the set of extreme points in its
domain.  For the set $\conv \Phi$, these extreme points can be realized by
deterministic quantizers based on likelihood
ratios~\cite{Tsitsiklis93}.  Consequently, we conclude that for
quantizers with binary outputs, the optimal cost is not decreased by
considering randomized quantizers. We conjecture that this statement
also holds beyond the binary case.
\end{enumerate}

Section~\ref{SecAsympBlock} is devoted to a more detailed study of
asymptotically optimal stationary quantizers.  In the meantime, we
turn to the question of whether stationary quantizers are optimal in
either finite-sample or asymptotic settings.

\comment{
\subsection{A simple illustration}

Consider two simple hypotheses $\Prob_0$ and $\Prob_1$ where $X$ takes
its values from the set $\{1,2\}$. In particular, $[\Prob_0(1)\; \Prob_0(2)] =
[0.8\; 0.2]$ and $[\Prob_1(1)\; \Prob_1(2) = [0.01\; 0.99]$.  The
priors are $\pi^0 = \pi^1 = 0.5$.
Table~\ref{table-example} shows that the approximation given 
by~\eqref{Eqn-approx} can provide a reasonable approximation to the 
cost for the optimal sequential test when $c$ is sufficiently small.

\begin{table}
\centerline{
\begin{tabular}{|r||r|r|r|r|r|r|r|r|r|r|}
\hline
c &0.01 & 0.009 & 0.008 & 0.007 & 0.006 & 0.005 & 0.004 & 0.003 & 0.002 
& 0.001 \\
\hline \hline
$J^*$ &  0.0320 &  0.0297 & 0.0269 & 0.0241 & 0.0212 & 0.0178 & 0.0145 
& 0.0112 & 0.0078 & 0.0042 \\
$\tilde{J}$ & 0.0302 &  0.0277 & 0.0250 & 0.0224 & 0.0196 & 0.0168 & 0.0139 
& 0.0108 & 0.0076 & 0.0041 \\
$\tilde{J}/J^*$ &  0.9447 & 0.9313 &  0.9313 & 0.9288 & 0.9268 & 0.9409 
& 0.9550 & 0.9682 & 0.9772 &  0.9737 \\
\hline
\end{tabular}}
\caption{
Comparison of the (exact) optimal cost $J^*$ computed by the dynamic 
programming method, and an approximation of the optimal cost
denoted by $\tilde{J}$ using Eq.~\eqref{Eqn-approx} as 
$c$ decreases.}
\label{table-example}
\end{table}
}

\section{Suboptimality of stationary designs}
\label{SecSuboptStat}

It was shown by Tsitsiklis~\cite{Tsitsiklis86} that optimal quantizers
$\phi_\timestep$ take the form of threshold rules based on the
likelihood ratio $f^1(X_\timestep)/f^0(X_\timestep)$.  Veeravalli et
al.~\cite{Veeravalli93,Veeravalli99} asked whether these rules can
always be taken to be stationary, a conjecture that has remained open.
In this section, we resolve this question with a negative answer in
both the finite-sample and asymptotic settings.

\subsection{Suboptimality in exact setting} 
\label{Subsec-exact}
We begin by providing a numerical counterexample for which stationary
designs are suboptimal.  Consider a problem in which $X \in
\mathcal{X} = \{1,2,3\}$ and the conditional distributions take the
form
\begin{equation*}
f^0(x) = \begin{bmatrix} \frac{8}{10} & \frac{1999}{10000} &
\frac{1}{10000}
\end{bmatrix} \; \mbox{and} \; f^1(x) = \begin{bmatrix} \frac{1}{3} &
\frac{1}{3} & \frac{1}{3} \end{bmatrix}.
\end{equation*}
Suppose that the prior probabilities are $\pi^1 = \frac{8}{100}$ and
$\pi^0 = \frac{92}{100}$, and that the cost for each sample is $c =
\frac{1}{100}$. 

If we restrict to binary quantizers (i.e., $\mathcal{U} = \{0,1\}$),
by the symmetric roles of the output alphabets 
there are only three possible deterministic quantizers:
\begin{enumerate}
\item Design A: $\phi_A(X_\timestep) = 0 \iff X_\timestep = 1$.  As a
result, the corresponding distribution for $U_\timestep$ is specified
by $f_{\phi_A}^0(u_n) = [\frac{4}{5} \; \frac{1}{5}]$ and
$f_{\phi_A}^1(u) = [\frac{1}{3}\; \frac{2}{3}]$.
\item Design B: ${\phi_B}(X_\timestep) = 0 \iff X_\timestep \in
\{1,2\}$.  The corresponding distribution for $U_\timestep$ is given
by $f_{\phi_B}^0(u) = [\frac{9999}{10000}\;\frac{1}{10000}]$ and
$f_{\phi_B}^1(u) = [\frac{2}{3}\; \frac{1}{3}]$.
\item Design C: $\phi_C(X_\timestep) = 0 \iff X_\timestep \in
\{1,3\}$.  The corresponding distribution for $U_\timestep$ is
specified by $f_{\phi_C}^0 \sim [\frac{8001}{10000}\;
\frac{1999}{10000}]$ and $f_{\phi_C}^1(u) = [\frac{2}{3} \;
\frac{1}{3}]$.
\end{enumerate}
Now consider the three stationary strategies, each of which uses only
one fixed design, A, B or C.  For any given stationary quantization
rule $\phi$, we have a classical centralized sequential problem, for
which the optimal cost (achieved by a sequential probability 
ratio test) can be computed using a dynamic-programming 
procedure~\cite{Waldetal48,Arrowetal49}.
Accordingly, for each stationary strategy, we compute the optimal cost
function $J$ for $10^6$ points on the $p$-axis by performing 300
updates of Bellman's equation (cf.~\cite{Bertsekas_dyn1}).  In all cases, the
difference in cost between the 299th and 300th updates is less than
$10^{-6}$. Let $J_A, J_B$ and $J_C$ denote the optimal cost function
for sequential tests using all A's, all B's, and all C's,
respectively.  When evaluated at $\pi^1=0.08$, these computations
yield $J_A = 0.0567, J_B = 0.0532$ and $J_C = 0.08$.
\comment{
\begin{table}
\begin{center}
\begin{tabular}{|c|c|c|c|c|}  
  \cline{1-5} & & & & \\
Method & $J_A(0.08)$ & $J_B(0.08)$ & $J_C(0.08)$ & $J_{*}(0.08)$ \\
& & & & \\
\cline{1-5} & & & & \\ Cost & $0.0567$ & $0.0532$ & $0.0800$ &
   $0.0528$ \\
    & & & & \\
    \cline{1-5}
\cline{1-5}
\end{tabular}
\end{center}
\caption{Numerically computed costs for the three stationary designs
$J_A$, $J_B$ and $J_C$, for the mixed design $J_*$.}
\label{TabExa}
\end{table}
}

Finally, we consider a non-stationary rule obtained by applying design
A for only the first sample, and applying design B for the remaining
samples. Again using Bellman's equation, we find that the cost for
this design is
\begin{eqnarray*}
J_* = \min \{\min\{\pi^1, 1-\pi^1 \}, c+ J_B(P(H=1|u_1 = 0))P(u_1 = 0) + && \\
J_B(P(H=1|u_1 = 1))P(u_1 = 1)\} = 0.052767,&&
\end{eqnarray*}
which is better than any of the stationary strategies. 

In this particular example, the cost $J^*$ of the non-stationary
quantizer yields a slim improvement (0.0004) over the best stationary
rule $J_B$.  This slim margin is due in part to the choice of a small
per-sample cost $c = 0.01$; however, larger values of $c$ do not yield
counterexample when using the particular distributions specified
above.  A more significant factor is that our non-stationary rule
differs from the optimal stationary rule $B$ only in its treatment of
the first sample.  This fact suggests that one might achieve better
cost by alternating between using design A and design B on the odd and
even samples, respectively.  Our analysis of the asymptotic setting in
the next section confirms this intuition.

\subsection{Asymptotic suboptimality for both deterministic
and randomized quantizers}

We now prove that in a broad class of examples, there is 
a range of prior probabilities for which stationary quantizer 
designs are suboptimal.  
Our result stems from the following observation: Lemma~\ref{PropApprox}
implies that in order
to achieve a small cost we need to choose a quantizer $\phi$ for which
the KL divergences $\asymKL{\phi} \defn D(f_\phi^0||f_\phi^1)$ 
and $\asymLK{\phi} \defn D(f_{\phi}^1||f_\phi^0)$ 
are both as large as possible. Due to the asymmetry of the 
KL divergence, however, these maxima are not necessarily
achieved by a single quantizer $\phi$.  This suggests that one
could improve upon stationary designs by applying different quantizers
to different samples, as the following lemma shows.
\begin{lemma}
\label{Prop-asymmetry}
Let $\phi_1$ and $\phi_2$ be any two quantizers.  If the following
inequalities hold
\begin{equation}
\label{EqnKeyAss}
\asymKL{\phi_1} < \asymKL{\phi_2} \mbox{ and } \asymLK{\phi_1} >
\asymLK{\phi_2}
\end{equation}
then there exists a non-empty interval $(U,V) \subseteq (0,+\infty)$
such that as $c\rightarrow 0$,
\begin{eqnarray*}
J_{\phi_1}^* \leq J_{\phi_1,\phi_2}^* \leq J_{\phi_2}^* &\mbox{if}&
\frac{\pi^0}{\pi^1} \leq U \\ 
J_{\phi_1,\phi_2}^* < \min \{ J_{\phi_1}^*, J_{\phi_2}^* \} -
\Theta(c\log c^{-1}) &\mbox{if}& \frac{\pi^0}{\pi^1} \in (U,V) \qquad
\\
J_{\phi_1}^* \geq J_{\phi_1,\phi_2}^* \geq J_{\phi_2}^* &\mbox{if}&
\frac{\pi^0}{\pi^1} \geq V,
\end{eqnarray*}
where $\Jfunc^*_{\phi_1,\phi_2}$
denotes the optimal cost of a sequential test that alternates between
using $\phi_1$ and $\phi_2$ on odd and even samples respectively.
\end{lemma}

\mybeginproof According to Lemma~\ref{PropApprox}, we have
\begin{eqnarray}
\label{Eqn-each}
\Jfunc_{\phi_i}^* & = & \biggr (\frac{\pi^0}{\asymKL{\phi_i}} +
\frac{\pi^1}{\asymLK{\phi_i}} \biggr ) c\log c^{-1} (1 + o(1)), \quad
i=0,1. \qquad
\end{eqnarray}
Now consider the sequential test that applies quantizers $\phi_1$ and
$\phi_2$ alternately to odd and even samples. Furthermore, let this test
consider two samples at a time. Let
$f_{\phi_1\phi_2}^0$ and $f_{\phi_1\phi_2}^1$ denote the induced
conditional probability distributions, jointly on the odd-even pairs
of quantized variables.  From the additivity of the KL divergence
and assumption~\eqref{EqnKeyAss}, there holds:
\begin{subequations}
\begin{eqnarray}
D(f_{\phi_1\phi_2}^0||f_{\phi_1\phi_2}^1) & = & \asymKL{\phi_1} +
\asymKL{\phi_2} > 2 \asymKL{\phi_1}\\
D(f_{\phi_1\phi_2}^1||f_{\phi_1\phi_2}^0) & = & \asymLK{\phi_1} +
\asymLK{\phi_2} < 2 \asymLK{\phi_1}.
\end{eqnarray}
\end{subequations}
Clearly, the cost of the proposed sequential test is an upper bound for
$\Jfunc_{\phi_1,\phi_2}^*$. Furthermore, the gap between this
upper bound and the true optimal cost is no more than $O(c)$.  Hence,
as in the proof of Lemma~\ref{PropApprox}, as $c\rightarrow 0$,
the optimal cost $\Jfunc_{\phi_1,\phi_2}^*$ can be written as
\begin{equation}
\label{Eqn-joint}
\biggr (\frac{2\pi^0}{\asymKL{\phi_1} + \asymKL{\phi_2}} +
\frac{2\pi^1}{\asymLK{\phi_1} + \asymLK{\phi_2}} \biggr ) c\log c^{-1}
(1+o(1)). \qquad
\end{equation}
From equations~\eqref{Eqn-each} and~\eqref{Eqn-joint}, simple
calculations yield the claim with
\begin{equation}
\label{Eqn-UV}
U = \frac{\asymKL{\phi_1}(\asymLK{\phi_1} -
\asymLK{\phi_2})(\asymKL{\phi_1}+\asymKL{\phi_2})}
{\asymLK{\phi_1}(\asymLK{\phi_1} + \asymLK{\phi_2}) (\asymKL{\phi_2} -
\asymKL{\phi_1})} < V = \frac{\asymKL{\phi_2}(\asymLK{\phi_1} -
\asymLK{\phi_2})(\asymKL{\phi_1}+\asymKL{\phi_2})}
{\asymLK{\phi_2}(\asymLK{\phi_1} + \asymLK{\phi_2}) (\asymKL{\phi_2} -
\asymKL{\phi_1})}.
\end{equation}
\myendproof

\noindent \myparagraph{Example:}
Let us return to the example provided in the previous
subsection.  Note that the two quantizers $\phi_A$ and $\phi_B$
satisfy assumption~\eqref{EqnKeyAss}, since
$D(f_{\phi_B}^0||f_{\phi_B}^1) = 0.4045 <
D(f_{\phi_A}^0||f_{\phi_A}^1) = 0.45 $ and
$D(f_{\phi_B}^1||f_{\phi_B}^0) = 2.4337 >
D(f_{\phi_A}^1||f_{\phi_A}^0) = 0.5108 $. Furthermore,  
both quantizers dominates $\phi_C$ in terms of KL divergences:
$D(f_{\phi_C}^0||f_{\phi_C}^1) = 0.0438$,
$D(f_{\phi_C}^0||f_{\phi_C}^1) = 0.0488$.
As a result, there exist a range of priors for which a sequential 
test using stationary quantizer design
(either $\phi_A$, $\phi_B$ or $\phi_C$ for all samples) is not
optimal. 

\begin{theorem} 
\label{ThmNewResult}
(a) Suppose that $\Phi$ is a finite collection of quantizers, and that
there is no single quantizer $\phi$ that dominates all other
quantizers in $\Phi$ in the sense that 
\begin{equation}
\asymKL{\phi} \; \geq \; \asymKL{\phi'} \quad \mbox{and} \quad
\asymLK{\phi} \; \geq \; \asymLK{\phi'} \qquad \mbox{for all} \quad
\phi' \in \Phi.
\end{equation}
Then there exists a non-empty range of prior probabilities for which
no stationary design based on a quantizer in $\Phi$ is optimal.

(b) For any non-deterministic $\phi$ in the randomized 
class $\conv \Phi$, there exists a non-stationary quantizer 
design that has strictly smaller sequential cost coefficient
than that of a stationary design based on $\phi$ for 
any choice of prior probabilities.
\end{theorem}
\begin{proof}
(a) Since there are a finite number of quantizers in $\Phi$ and no quantizer
dominates all others, the interval $(0,\infty)$ is divided into at
least two adjacent non-empty intervals, each of which corresponds to a
range of prior probability ratios $\pi^0/\pi^1$ for which a quantizer
is strictly optimal (asymptotically) among all stationary designs.
Let them be $(\delta_1,\delta)$ and $(\delta, \delta_2)$, for two
quantizers, namely, $\phi_1$ and $\phi_2$. In particular, $\delta$ is
the value for $\pi^0/\pi^1$ for which the sequential cost coefficients
are equal---viz. $G_{\phi_1} = G_{\phi_2}$---which happens only if
assumption~\eqref{EqnKeyAss} holds.  Some calculations verify that
\begin{equation}
\label{EqnDefnDelta}
\delta = \frac{\asymKL{\phi_1}\asymKL{\phi_2}
(\asymLK{\phi_2}-\asymLK{\phi_1})} {\asymLK{\phi_1}\asymLK{\phi_2}
(\asymKL{\phi_1}-\asymKL{\phi_2})}.
\end{equation}
By Lemma~\ref{Prop-asymmetry}, a non-stationary design obtained by alternating
between $\phi_1$ and $\phi_2$ has smaller sequential cost than both
$\phi_1$ and $\phi_2$ for $\pi^0/\pi^1 \in (U,V)$, where $U$ and $V$
are given in equation~\eqref{Eqn-UV}. Since it can be verified that
$\delta$ as defined~\eqref{EqnDefnDelta} belongs to the interval
$(U,V)$, we conclude that for $\pi^0/\pi^1 \in (U,V)\cap
(\delta_1,\delta_2)$, this non-stationary design has smaller cost than
any stationary design using $\phi \in \Phi$.

(b) Let $\phi \in \conv \Phi$ be a randomized quantizer (i.e., at each
step choose with non-zero probabilities $w_1,\ldots, w_k$ from
quantizers $\phi_1,\ldots,\phi_k \in \Phi$, respectively, where
$\sum_{i=1}^{k}w_i = 1$). Clearly, the density induced by $\phi$
satisfy: $f_\phi^0 = \sum_{i=1}^{k}w_i f_{\phi_i}^0$ and $f_\phi^1 =
\sum_{i=1}^{k}w_i f_{\phi_i}^1$. Due to strict convexity of the KL
divergence functional with respect jointly to the two density
arguments~\cite{Cover91}, by Jensen's inequality we have:
$\mykull_\phi^0 < \sum_{i=1}^{k}w_i\mykull_{\phi_i}^0$ and
$\mykull_\phi^1 < \sum_{i=1}^{k}w_i\mykull_{\phi_i}^1$.  Since
$\mykull_{\phi_i}^0$ and $\mykull_{\phi_i}^1$ are bounded from above
uniformly for all $\phi_i\in \Phi$, it is possible to approximate
$(w_1,\ldots,w_k)$ by rational numbers of the form $(q_1/N,
q_2/N,\ldots, q_k/N)$ for some natural numbers $q_1,\ldots, q_k$ and
$N$ satisfying $\sum_{i=1}^{k} q_i = N$ such that
\begin{eqnarray*}
\mykull_\phi^0 & < & \sum_{i=1}^{k}q_i \mykull_{\phi_i}^0/N \\
\mykull_\phi^1 & < & \sum_{i=1}^{k}q_i \mykull_{\phi_i}^1/N.
\end{eqnarray*}
Now consider the non-stationary quantizer that applies $\phi_1$ for
$q_1$ steps, then $\phi_2$ for $q_2$ steps and so on, up to $\phi_k$
for $q_k$ steps, yielding a total of $N$ steps, and then repeats this
sequence starting again at step $N+1$.  By construction, this
non-stationary quantizer has a smaller cost than that of quantizer
$\phi$ for any choice of prior.
%
\end{proof}
\noindent \myparagraph{Remarks:} (i) It is worth emphasizing
the assumption that the class $\Phi$ is finite is crucial in 
part a) of the theorem. We do not know if this result can
be extended to the case in which $\Phi$ is infinite. 
(ii) Part b) shows that
any stationary randomized quantizer is always dominated by some non-stationary
one.  Actually, a stronger result can be proved at least
for binary quantizers (see Cor.~\ref{Cor-binary}): 
for any given choice of prior probability, any stationary randomized 
quantizer is dominated by a stationary deterministic quantizer.
(iii) It is interesting to contrast the Bayesian formulation of the problem
of quantizer design with the Neyman-Pearson formulation.  Our results
on the suboptimality of stationary quantizer design in the Bayesian
formulation repose on the asymmetry of the Kullback-Leibler
divergence, as well as the sensitivity of the optimal quantizers on
the prior probability.  We note that Mei~\cite{Mei-thesis} (see p. 58)
considered the Neyman-Pearson formulation of this problem.  In this
formulation, it can be shown that for all sequential tests for which
the Type 1 and Type 2 errors are bounded by $\alpha$ and $\beta$,
respectively, then as $\alpha + \beta \rightarrow 0$, the expected
stopping time $\E_0 N$ under hypothesis $H=0$ is asymptotically
minimized by applying a stationary quantizer $\phi^*$ that maximizes
$D(f_\phi^0 ||f_\phi^1)$.  Similarly, the expected stopping time $\E_1
N$ under hypothesis $H=1$ is asymptotically minimized by the
stationary quantizer $\phi^{**}$ that maximizes
$D(f_\phi^1||f_\phi^0)$~\cite{Mei-thesis}.  In this context, the
example in subsection~\ref{Subsec-exact} provides a case in which the
asymptotically minimal KL divergences $\phi^*$ and $\phi^{**}$ are not
the same, due to the asymmetry, which suggests that there may not
exist a stationary quantizer that simultaneously minimizes both $\E_1
N$ and $\E_0 N$.

\subsection{Asymptotic suboptimality in multiple sensor setting}
Our analysis thus far has established that with a single sensor per
time step ($\dims = 1$), applying multiple quantizers to different
samples can reduce the sequential cost.  As pointed out by one of
the referees, it is natural to ask whether the same phenomenon
persists in the case of multiple sensors ($\dims > 1$).  In this
section, we show that the phenomenon does indeed carry over, more
specifically by providing an example in which stationary strategies
are still sub-optimal in comparison to non-stationary ones.  The key
insight is that we have only a fixed number of dimensions, whereas as
$c\rightarrow 0$ we are allowed to take more samples, and each sample
can act as an extra dimension, providing more flexibility for
non-stationary strategies.

Suppose that the observation vector $X_n$ at time $n$ is
$\dims$-dimensional, with each component corresponding to a sensor in
a typical decentralized setting.  Suppose that the observations from
each sensor are assumed to be independent and identically distributed
according to the conditional distributions defined in our earlier
example (see Section~\ref{Subsec-exact}).  Of interest are the optimal
deterministic binary quantizer designs for all $\dims$
sensors. Although there are three possible choices $\phi_A$, $\phi_B$
and $\phi_C$ for each sensor, the quantizer $\phi_C$ is dominated by
the other two, so each sensor should choose either $\phi_A$ and
$\phi_B$. Suppose that among these sensors, a subset of size $\ind$
choose $\phi_A$ and whereas the remaining $\dims - \ind$ sensors
choose $\phi_B$ for $0\leq k \leq \dims$.  We thus have $\dims + 1$
possible stationary designs to consider.  For each $\ind$, the
sequential cost coefficient corresponding to the associated stationary
design takes the form
\begin{equation}
G_k \defn \frac{\pi^0}{\ind\asymKL{\phi_A} + (\dims -
\ind)\asymKL{\phi_B}} + \frac{\pi^1}{\ind\asymLK{\phi_A} + (\dims -
\ind)\asymLK{\phi_B}}.
\end{equation}

Now consider the following non-stationary design: the first sensor
alternates between decision rules $\phi_A$ and $\phi_B$, while the
remaining $\dims - 1$ sensors simply apply the stationary design based
on $\phi_B$.  For this design, the associated sequential cost
coefficient is given by
\begin{equation}
\label{EqnNewSequentialCostCoeff}
G \defn \frac{2\pi^0}{\asymKL{\phi_A} + (2\dims - 1)\asymKL{\phi_B}} +
\frac{2\pi^1}{\asymLK{\phi_A} + (2\dims - 1)\asymLK{\phi_B}}.
\end{equation}
Consider the interval $(U,V)$, where the interval has endpoints
\[U = \frac{\asymLK{\phi_B} - \asymLK{\phi_A}}{\asymKL{\phi_A}-
\asymKL{\phi_B}} \frac{\asymKL{\phi_A} + (2\dims-1)\asymKL{\phi_B}}
{\asymLK{\phi_A} + (2\dims-1)\asymLK{\phi_B}}
\frac{\asymKL{\phi_B}}{\asymLK{\phi_B}} < V = \frac{\asymLK{\phi_B} -
\asymLK{\phi_A}}{\asymKL{\phi_A}- \asymKL{\phi_B}}
\frac{\asymKL{\phi_A} + (2\dims-1)\asymKL{\phi_B}} {\asymLK{\phi_A} +
(2\dims-1)\asymLK{\phi_B}} \frac{\asymKL{\phi_A} + (\dims -
1)\asymKL{\phi_B}} {\asymLK{\phi_A} + (\dims - 1)\asymLK{\phi_B}}.
\]
Straightforward calculations yield that for any prior likelihood
$\pi^0/\pi^1 \in (U,V)$, the minimal cost over stationary designs
$\min_{\ind = 0,\ldots,\dims} G_k$ is strictly larger than the
sequential cost $G$ of the non-stationary design, previously defined
in equation~\eqref{EqnNewSequentialCostCoeff}.

\section{On asymptotically optimal blockwise stationary designs}
\label{SecAsympBlock}

Despite the possible loss in optimality, it is useful to consider 
some form of stationarity in order to reduce computational complexity 
of the optimization and decision process.  In this section, we consider
the class of \emph{blockwise stationary} designs, meaning that there
exists some natural number $T$ such that $\phi_{T+1} = \phi_1,
\phi_{T+2} = \phi_2,$ and so on.  For each $T$, let $C_T$ denote the
class of all blockwise stationary designs with period $T$.  We assume
throughout the analysis that each decision rule $\phi_n
\;(n=1,\ldots,T)$ satisfies conditions~\eqref{Eqn-Ass1} and
~\eqref{Eqn-Ass2}.
Thus, as $T$ increases, we have a hierarchy of increasingly rich
quantizer classes that will be seen to yield progressively better 
approximations to the optimal solution.

For a fixed prior $(\pi^0,\pi^1)$ and $T>0$, let 
$(\phi_1,\ldots,\phi_T)$ denote a quantizer design 
in $C_T$. As before, the cost $J_{\phi}^*$ of an asymptotically optimal
sequential test using this quantizer design is
of order $c\log c^{-1}$ with the sequential
cost coefficient
\begin{equation}
\label{Eqn-stationary}
 G_\phi = \frac{T\pi^0}{\asymKL{\phi_1} +\ldots+ \asymKL{\phi_T}} +
\frac{T\pi^1}{\asymLK{\phi_1} +\ldots+ \asymLK{\phi_T}}.
\end{equation}
$G_\phi$ is a function of the vector of probabilities introduced
by the quantizer: $(f_\phi^0(.),f_\phi^1(.))$.
We are interested in the properties of a quantization 
rule $\phi$ that minimizes $J_\phi^*$.

\comment{
We begin with a simple result on the structure of asymptotically
optimal quantizer designs:
\begin{proposition}
For a fixed prior $(\pi^0,\pi^1)$ and $T>0$, let
$(\phi_1^*,\ldots,\phi_T^*)$ be the optimal quantizer design among
those in $C_T$.
\begin{enumerate}
\item[(a)] For any $1\leq i,j \leq T$, there holds $(\asymKL{\phi_i^*}
  - \asymKL{\phi_j^*})(\asymLK{\phi_i^*} - \asymLK{\phi_j^*}) \leq 0.$
\item[(b)] If there is a pair $(i,j)$ such that $(\asymKL{\phi_i^*} -
\asymKL{\phi_j^*})(\asymLK{\phi_i^*} - \asymLK{\phi_j^*}) < 0$ then
there exists a prior $\pi^0$ so that $(\phi_1^*,\ldots,\phi_T^*)$ is
not optimal in $C_T$.
\end{enumerate}
\end{proposition}
\begin{proof}
To establish claim (a), suppose that for some pair $(i,j)$, there
holds $(\asymKL{\phi_i^*} - \asymKL{\phi_j^*})(\asymLK{\phi_i^*} -
\asymLK{\phi_j^*}) > 0.$ Without loss of generality, assume that
$\asymKL{\phi_i^*} > \asymKL{\phi_j^*}$ and $\asymLK{\phi_i^*} >
\asymLK{\phi_j^*})$.  The asymptotic cost~\eqref{Eqn-stationary}
can then be improved by replacing quantizer design $\phi_j^*$ at the
periodic index $j$ by the design $\phi_i^*$.  The proof of
part (b) is similar to that of Proposition~\ref{Prop-asymmetry}.
\end{proof}
}

It is well known that there exist optimal
quantizers---\emph{when unrestricted}--- that can be expressed as threshold
rules based on the log likelihood ratio (LLR)~\cite{Tsitsiklis86}. 
Our counterexamples in
the previous sections imply that the thresholds need not be stationary
(i.e., the threshold may differ from sample to sample).  In the
remainder of this section, we addresses a partial converse to this
issue: specifically, if we restrict ourselves to stationary (or
blockwise stationary) quantizer designs, then there exists an optimal
design consisting of LLR-based threshold rules.

It turns out that the analysis for the case $T > 1$ can be
reduced to an analysis that is closely related to our earlier 
analysis for $T = 1$.  Indeed, consider the sequential cost 
coefficient for the time step $\timestep = 1$, where the rules 
for the other time steps are held fixed.  From~\eqref{Eqn-stationary}
we have
\[G_\phi = \frac{T\pi^0}{\asymKL{\phi_1} + s_0} +
\frac{T\pi^1}{\asymLK{\phi_1} + s_1},\]
for non-negative constants $s_0$ and $s_1$.  
As we will show, our earlier analysis of the sequential cost
coefficient, in which $s_0 = s_1 = 0$, carries through to 
the case in which these values are non-zero.  This allows us 
to provide (in Theorem~\ref{ThmPartialConv}) a characterization 
of the optimal blockwise stationary quantizer. 

\begin{definition} The quantizer design 
function $\phi: \Xspace \rightarrow \Uspace$ is said to be a 
\emph{likelihood ratio threshold
rule} if there are thresholds $d_0 = -\infty < d_1 <\ldots < d_{\Umax}
= +\infty$, and a permutation $(u_1,\ldots,u_{\Umax})$ of
$(0,1,\ldots,\Umax-1)$ such that for $l = 1,\ldots,\Umax$, with
$\Prob_0$-probability 1, we have:
\[\phi(X) = u_l \;\mbox{if}\;\; d_{l-1} 
\leq f^1(X)/f^0(X) \leq d_{l},\]
When $f^1(X)/f^0(X) = d_{l-1}$, set $\phi(X) =
u_{l-1}$ or $\phi(X) = u_{l}$ with $\Prob_0$-probability
1.\footnote{This last requirement of the definition is termed
the \emph{canonical} likelihood ratio quantizer by
Tsitsiklis~\cite{Tsitsiklis93-extremal}. Although one could consider
performing additional randomization when there are ties, our later results
(in particular, Lemma~\ref{Lem-Quasiconcave}) establish that in this
case, randomization will not further decrease the optimal cost
$J_\phi^*$.}
\end{definition}

Previous work on the extremal properties of likelihood ratio based
quantizers guarantees that the Kullback-Leibler divergence is
maximized by a LLR-based quantizer~\cite{Tsitsiklis93-extremal}.
In our case, however, the sequential cost coefficient $G_\phi$
involves a pair of KL divergences,
$\asymKL{\phi}$ and $\asymLK{\phi}$, which are related to one
another in a nontrivial manner. Hence, establishing asymptotic
optimality of LLR-based rules for this cost function does not follow
from existing results, but rather requires further understanding of
the interplay between these two KL divergences.

The following lemma concerns certain ``unnormalized'' variants of the
Kullback-Leibler (KL) divergence.  Given vectors $a =(a_0, a_1)$ and
$b = (b_0, b_1)$, we define functions $\extKLplain$ and $\extLKplain$
mapping from $\real_+^4$ to the real line as follows:
\begin{subequations}
\begin{eqnarray}
\extKL{a}{b} & \defn & a_0\log\frac{a_0}{a_1} + b_0\log\frac{b_0}{b_1}
\\ \extLK{a}{b} & \defn & a_1\log\frac{a_1}{a_0} +
b_1\log\frac{b_1}{b_0}.
\end{eqnarray}
\end{subequations}
These functions are related to the standard (normalized) KL divergence
via the relations $\extKL{a}{1-a} \equiv D(a_0,a_1)$, and
\mbox{$\extLK{a}{1-a} \equiv D(a_1,a_0)$.}

\begin{lemma}
\label{Lem-Ineq} 
For any positive scalars $a_1, b_1, c_1, a_0, b_0, c_0$ such that
$\frac{a_1}{a_0} < \frac{b_1}{b_0} < \frac{c_1}{c_0}$, \emph{at least
one} of the two following conditions must hold:
\begin{subequations}
\begin{eqnarray}
\extKL{a}{b+c} > \extKL{b}{c+a} & \mbox{and} & \extLK{a}{b+c} >
\extKL{b}{c+a}, \;\; \mbox{or}  \\
\extKL{c}{a+b} > \extKL{b}{c+a} & \mbox{and} & \extLK{c}{a+b} >
\extKL{b}{c+a}.
\end{eqnarray}
\end{subequations}
\end{lemma}
This lemma implies that under certain conditions on the
ordering of the probability ratios, one can increase \emph{both} 
KL divergences by re-quantizing. This insight is used in the following 
lemma to establish that the optimal quantizer $\phi$ behaves almost 
like a likelihood ratio rule.  To state the result,
recall that the \emph{essential supremum} is the infimum of the 
set of all $\eta$ such that $f(x) \leq \eta$ for $\Prob_0$-almost all 
$x$ in the domain, for a measurable function $f$.
\begin{lemma}
\label{Prop-PseudoLR}
If $\phi$ is an asymptotically optimal quantizer, then for all pairs
$(u_1,u_2) \in \Uspace$, $u_1\neq u_2$, there holds:
\[
\frac{f^1(u_1)}{f^0(u_1)} \notin \biggr (\ess \inf_{x: \phi(x) = u_2}
\frac{f^1(x)}{f^0(x)}, \; \ess \sup_{x:\phi(x) = u_2}
\frac{f^1(x)}{f^0(x)} \biggr).
\]
\end{lemma}
Note that a likelihood ratio rule guarantees something stronger: For
$\Prob_0$-almost all $x$ such that $\phi(x) = u_1$, $f^1(x)/f^0(x)$
takes a value either to the left or to the right, but not to both sides,
of the interval specified above. 
\comment{
As we shall show,
the proof that there exists an optimal LLR-based rule
turns out to reduce to the problem of showing that the sequential
cost coefficient
$G_\phi$ is a \emph{quasiconcave} function with respect
to $(f_\phi^0(.), f_\phi^1(.))$. 
Since the minima of a quasiconcave function are extreme points
of the function's domain~\cite{Boyd-book}, and the extreme points in the 
quantizer space are LLR-based rules (Prop. 3.2 
of~\cite{Tsitsiklis93-extremal}), we 
deduce that there exists an optimal quantizer that is 
LLR-based.}

Lemma~\ref{Lem-Quasiconcave} stated below essentially guarantees 
quasiconcavity of $G_\phi$ for the case of binary quantizers.  
To state the result, let $F:[0,1]^2
\rightarrow R$ be given by
\begin{equation}
\label{EqnDefnF}
F(a_0,a_1) = \frac{c_0}{D(a_0,a_1) + d_0} +  \frac{c_1}{D(a_1,a_0) +
d_1}.
\end{equation}
\begin{lemma}
\label{Lem-Quasiconcave}
For any non-negative constants $c_0,c_1,d_0,d_1$, the function $F$
defined in~\eqref{EqnDefnF} is quasiconcave.
\end{lemma}
We provide a proof of this result in the Appendix. An immediate
consequence of Lemma~\ref{Lem-Quasiconcave} is that LLR-based quantizers
exist for the class of randomized quantizers with binary outputs. 
\begin{cor}
\label{Cor-binary}
Restricting to the class of (blockwise) stationary binary quantizers,
there exists an asymptotically optimal quantizer 
$\phi$ that is a (deterministic) likelihood ratio threshold rule.
\end{cor}
\mybeginproof
Let $\phi$ is a (randomized) binary quantizer. The sequential cost
coefficient can be written as $G_\phi =
F(f_\phi^0(0),f_\phi^1(0))$. The set of $\{(f_\phi^0(0),f_\phi^1(0)\}$
for all $\phi$ is a convex set whose extreme points can be realized
by deterministic likelihood ratio threshold rules 
(Prop. 3.2 of~\cite{Tsitsiklis93-extremal}).  
Since the minimum of a quasiconcave function must lie at one 
such extreme point~\cite{Boyd-book}, the corollary is immediate 
as a consequence of Lemma~\ref{Lem-Quasiconcave}.
\myendproof

It turns out that the same statement can also be proved for deterministic
quantizers with arbitrary output alphabets:
\begin{theorem}
\label{ThmPartialConv}
Restricting to the class of (blockwise) stationary and deterministic
decision rules, then there exists an asymptotically optimal quantizer
$\phi$ that is a likelihood ratio threshold rule.
\end{theorem}
We present the full proof of this theorem in
the Appendix.  The proof exploits both
Lemma~\ref{Prop-PseudoLR} and Lemma~\ref{Lem-Quasiconcave}.

\section{Discussion}
\label{SecDiscussion}
In this paper, we have studied the problem of sequential decentralized
detection.  For quantizers with neither local memory nor feedback
(Case A in the taxonomy of Veeravalli et al.~\cite{Veeravalli93}), we
have established that stationary designs need not be optimal in
general.  Moreover, we have shown that in the asymptotic setting
(i.e., when the cost per sample goes to zero), there is a class of
problems for which there exists a range of prior probabilities over
which stationary strategies are suboptimal.

There are a number of open questions raised by the analysis in this
paper.  First, our analysis has established only that the best
stationary rule chosen from a finite set of deterministic quantizers
need not be optimal.  Is there a corresponding example with an
infinite number of deterministic stationary quantizer designs for
which none is optimal?  Second, Corollary~\ref{Cor-binary} establishes
the optimality of likelihood ratio rules for randomized decision rules
that produce binary outputs.  This proof was based on the
quasiconcavity of the function $G_\phi$ that specifies the asymptotic
sequential cost coefficient. Is this function $G_\phi$ also
quasiconcave for quantizers other than binary ones?  Such
quasiconcavity would extend the validity of
Theorem~\ref{ThmPartialConv} for the general class of randomized
quantizers.

\subsection*{Acknowledgements}
The research described in this paper was supported in part by National
Science Foundations Grants 0412995 (MJ) and DMS-0605165 (MW).  We thank
the reviewers for their comments and suggestions which helped to improve
the manuscript.

\appendix

\mysection{Proof of Lemma~\ref{Lem-Ineq}}
\label{AppLem_Ineq}

By renormalizing, we can assume w.l.o.g. that $a_1+b_1+c_1 = a_0 + b_0
+ c_0 = 1$.  Also w.l.o.g, assume that $b_1 \geq b_0$. Thus, $c_1 >
c_0$ and $a_1 < a_0$. Replacing $c_1 = 1 - a_1 - b_1$ and $c_0 = 1 -
a_0 - b_0$, the inequality $c_1/c_0>b_1/b_0$ is equivalent to $a_1 <
a_0 b_1/b_0 - (b_1 - b_0)/b_0$.

\begin{figure}
\begin{center}
\psfrag{#a0#}{$a_0$}
\psfrag{#a1#}{$a_1$}
\psfrag{#up0#}{$1-b_0$}
\psfrag{#up1#}{$1-b_1$}
\psfrag{#down0#}{$(b_1-b_0)/b_1$}
\psfrag{#A#}{$A$}
\widgraph{.37\textwidth}{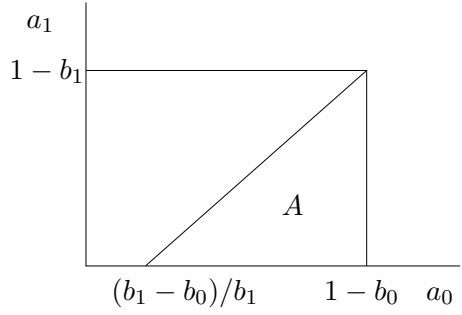}
\end{center}
\caption{Illustration of the domain $A$.}
\label{Fig-A}
\end{figure}

We fix values of $b$, and consider varying $a \in A$, where $A$
denotes the domain for $(a_0,a_1)$ governed by the following
equality and inequality constraints:
$0 < a_1 < 1 - b_1; \;
 0 < a_0 < 1 - b_0; \;
 a_1 < a_0$ and
\begin{equation}
\label{boundary-1}
a_1 < a_0 b_1/b_0 - (b_1 - b_0)/b_0.
\end{equation}
Note that the third constraint ($a_1 < a_0$) is redundant due to the
other three constraints. In particular, constraint~\eqref{boundary-1}
corresponds to a line passing through $((b_1-b_0)/b_1, 0)$ and
$(1-b_0,1-b_1)$ in the $(a_0,a_1)$ coordinates. As a result, $A$ is
the interior of the triangle defined by this line and two other lines
given by $a_1 = 0$ and $a_0 = 1-b_0$ (see Figure~\ref{Fig-A}).

Since both $\extKL{a}{1-a}$ and $\extLK{a}{1-a}$ correspond to KL
divergences, they are convex functions with respect to $(a_0,a_1)$.
In addition, the derivatives with respect to $a_1$ are
$\frac{a_1-a_0}{a_1(1-a_1)} < 0$ and
$\log\frac{a_1(1-a_0)}{a_0(1-a_1)} < 0$, respectively.  Hence, both
functions can be (strictly) bounded from below by increasing $a_1$
while keeping $a_0$ unchanged, i.e., by replacing $a_1$ by $a'_1$ so
that $(a_0,a'_1)$ lies on the line given by~\eqref{boundary-1}, which
is equivalent to the constraint $c_1/c_0 = b_1/b_0$.  Let $c'_1 = 1 -
b_1 - a'_1$; then $c'_1/c_0 = b_1/b_0$. Our argument has established
inequalities (a) and (b) in the following chain of inequalities:
\begin{subequations}
\begin{eqnarray}
\extLK{a}{b+c} & \stackrel{(a)}{>} & a'_1 \log \frac{a'_1}{a_0} + (b_1
+ c'_1)\log \frac{b_1 + c'_1}{b_0 + c_0} \\
& \stackrel{(b)}{=} & a'_1 \log \frac{a'_1}{a_0} +
c'_1\log\frac{c'_1}{c_0} + b_1\log\frac{b_1}{b_0} \\
 & \stackrel{(c)}{\geq} & (a'_1+c'_1)\log\frac{a'_1+c'_1}{a_0+c_0} +
b_1\log\frac{b_1}{b_0} \\ & = & \extLK{a+c}{b},
\end{eqnarray}
\end{subequations}
inequality (c) follows from an application of the log-sum
inequality~\cite{Cover91}.  A similar conclusion holds for
$\extKL{a}{b+c}$.

%
\mysection{Proof of Lemma~\ref{Prop-PseudoLR}}

Suppose the opposite is true, that there exist two sets $S_1$, $S_2$
with positive $\Prob_0$-measure such that $\phi(X) = u_2$ for any $X
\in S_1 \cup S_2$, and
\begin{equation}
\frac{f^1(S_1)}{f^0(S_1)} < \frac{f^1(u_1)}{f^0(u_1)} <
\frac{f^1(S_2)}{f^0(S_2)}.
\end{equation}
By reassigning $S_1$ or $S_2$ to the quantile $u_1$, we are guaranteed
to have a new quantizer $\phi'$ such that $\asymKL{\phi'} >
\asymKL{\phi^*}$ and $\asymLK{\phi'} > \asymLK{\phi^*}$, thanks
to Lemma~\ref{Lem-Ineq}.  As a result, $\phi'$ has a smaller
sequential cost $J_{\phi'}^*$, which is a contradiction.


\mysection{Proof of Lemma~\ref{Lem-Quasiconcave}}
\label{Proof-quasiconcave}

The proof of this lemma is conceptually straightforward, but
the algebra is involved.
To simplify the notation, we replace $a_0$ by $x$, $a_1$ by $y$, the
function $D(a_0,a_1)$ by $f(x,y)$, and the function $D(a_1,a_0)$ by
$g(x,y)$.  Finally, we
assume that $d_0= d_1 = 0$; the proof will reveal that this case is
sufficient to establish the more general result with arbitrary non-negative
scalars $d_0$ and $d_1$. 

We have $f(x,y) = x\log (x/y) + (1-x)\log [(1-x)/(1-y)]$ and $g(x,y) =
y\log (y/x) + (1-y)\log [(1-y)/(1-x)]$.  Note that both $f$ and $g$
are convex functions and are non-negative in their domains, and
moreover that we have $F(x,y) = c_0/f(x,y) + c_1/g(x,y)$.  In order to
establish the quasiconcavity of $F$, it suffices to show that for any
$(x,y)$ in the domain of $F$, for any vector $h = [h_0 \; h_1] \in
\real^2$ such that $h^T \nabla F(x, y) = 0$, there holds
\begin{equation}
\label{Eqn-QuasiCond}
h^T \nabla^2 F(x,y) \; h < 0
\end{equation} 
(see Boyd and Vandenberghe~\cite{Boyd-book}).  Here we adopt the
standard notation of $\nabla F$ for the gradient vector of $F$, and
$\nabla^2 F$ for its Hessian matrix.  We also use $F_x$ to denote the
partial derivative with respect to variable $x$, $F_{xy}$ to denote
the partial derivative with respect to $x$ and $y$, and so on.

We have $\nabla F = -\frac{c_0 \nabla f}{f^2} - \frac{c_1 \nabla
g}{g^2}$.  Thus, it suffices to prove relation~\eqref{Eqn-QuasiCond}
for vectors of the form
\begin{equation*}
h = \begin{bmatrix} \left(-\frac{c_0 f_y}{f^2} - \frac{c_1
g_y}{g^2}\right) & \left(\frac{c_0 f_x}{f^2} + \frac{c_1 g_x}{g^2}
\right)
\end{bmatrix}^T.
\end{equation*}
It is convenient to write $h = c_0 v_0 + c_1 v_1$, where $v_0 =
[-f_y/f^2 \;\;\; f_x/f^2]^T$ and $v_1 = [-g_y/g^2 \;\;\; g_x/g^2]^T$.

The Hessian matrix $\nabla^2 F$ can be written as $\nabla^2 F = c_0
H_0+ c_0 H_1$, where
\[ H_0 = -\frac{1}{f^3} \biggr [ \begin{array}{cc}
\derivtwo{f}{x}{x}f - 2\derivone{f}{x}^2 & \derivtwo{f}{x}{y}f -
2\derivone{f}{x}\derivone{f}{y}\\ \derivtwo{f}{x}{y}f -
2\derivone{f}{x}\derivone{f}{y} & \derivtwo{f}{y}{y}f -
2\derivone{f}{y}^2 \end{array} \biggr ],\]
and 
\[ H_1 = -\frac{1}{g^3} \biggr [ \begin{array}{cc}
\derivtwo{g}{x}{x}g - 2\derivone{g}{x}^2 & \derivtwo{g}{x}{y}g - 
2\derivone{g}{x}\derivone{g}{y}\\
\derivtwo{g}{x}{y}g - 2\derivone{g}{x}\derivone{g}{y} & \derivtwo{g}{y}{y}g - 
2\derivone{g}{y}^2 \end{array} \biggr ].\]

Now observe that
\begin{eqnarray*}
h^T \nabla^2 F h & = & (c_0 v_0 + c_1 v_1)^T(c_0 H_0 + c_1 H_1) (c_0
v_0 + c_1 v_1),
\end{eqnarray*}
which can be simplified to
\begin{multline*}
h^T \nabla^2 F h = c_0^3 v_0^T H_0 v_0 + c_1^3 v_1^T H_1 v_1 + c_0^2
c_1(2v_0^T H_0 v_1 + v_0^T H_1 v_0) + c_0 c_1^2(2v_0^T H_1 v_1 + v_1^T
H_0 v_1).
\end{multline*}
This function is a polynomial in $c_0$ and $c_1$, which are restricted
to be non-negative scalars (at least one of which is assumed to be non-zero).  
Therefore, it suffices to prove that all
the coefficients of this polynomial (with respect to $c_0$ and $c_1$)
are strictly negative. In particular, we shall show that
\begin{enumerate}
\item[(i)] $v_0^T H_0 v_0 \leq 0$, and
\item[(ii)] $2v_0^T H_0 v_1 + v_0^T H_1 v_0 \leq 0$,. 
\end{enumerate}
where in both cases equality occurs only if $x = y$, which is 
outside of the domain of $F$. The strict negativity of the other 
two coefficients follows from entirely analogous arguments.

First, some straightforward algebra shows that inequality (i) is
equivalent to the relation
\begin{equation*}
\derivtwo{f}{x}{x} \derivone{f}{y}^2 +
\derivtwo{f}{y}{y}\derivone{f}{x}^2 \geq
2\derivone{f}{x}\derivone{f}{y}\derivtwo{f}{x}{y}.
\end{equation*}
But note that $f$ is a convex function, so
$\derivtwo{f}{x}{x}\derivtwo{f}{y}{y}\geq \derivtwo{f}{x}{y}^2$.  Hence,
we have
\begin{equation*}
\derivtwo{f}{x}{x} \derivone{f}{y}^2 +
\derivtwo{f}{y}{y}\derivone{f}{x}^2 \; \stackrel{(a)}{\geq} \;
2\sqrt{\derivtwo{f}{x}{x}\derivtwo{f}{y}{y}}|\derivone{f}{x}\derivone{f}{y}|
\; \stackrel{(b)}{\geq} \;
2\derivone{f}{x}\derivone{f}{y}\derivtwo{f}{x}{y},
\end{equation*}
thereby proving (i).  (In this argument, inequality (a) follows from
the fact that $a^2 + b^2 \geq 2 a b$, whereas inequality (b) follows
from the strict convexity of $f$. Equality occurs
only if $x = y$.)

Regarding (ii), some further algebra reduces it to the inequality
\begin{equation}
\label{Eqn-Algebra-1}
G_1 + G_2 - G_3 \geq 0,
\end{equation}
where
\begin{eqnarray*}
G_1 & = & 2(\derivone{f}{y}\derivone{g}{y}\derivtwo{f}{x}{x} +
  \derivone{f}{x}\derivone{g}{x}\derivtwo{f}{y}{y} -
  (\derivone{f}{y}\derivone{g}{x} + \derivone{f}{x}\derivone{g}{y})
  \derivtwo{f}{x}{y}), \\ G_2 & = & \derivone{f}{y}^2
  \derivtwo{g}{x}{x} + \derivone{f}{x}^2 \derivtwo{g}{y}{y} -
  2\derivone{f}{x}\derivone{f}{y}\derivtwo{g}{x}{y}, \\ G_3 & = &
  \frac{2}{g}(\derivone{f}{y}\derivone{g}{x} -
  \derivone{f}{x}\derivone{g}{y})^2.
\end{eqnarray*}

At this point in the proof, we need to exploit specific information
about the functions $f$ and $g$, which are defined in terms of KL
divergences.  To simplify notation, we let $u = x/y$ and
$v=(1-x)/(1-y)$.  Computing derivatives, we have
\begin{eqnarray*}
\derivone{f}{x}(x,y) & = & \log (x/y) - \log((1-x)/(1-y)) = \log
(u/v), \\
\derivone{f}{y}(x,y) & = & (1-x)/(1-y) - x/y = v - u, \\
\derivone{g}{x}(x,y) & = & (1-y)/(1-x) - y/x = 1/v - 1/u, \\
\derivone{g}{y}(x,y) & = & \log (y/x) - \log((1-y)/(1-x)) = \log
(v/u), \\
\nabla^2 f(x,y) & = & \biggr [ \begin{array}{cc} \frac{1}{x(1-x)} &
    -\frac{1}{y(1-y)} \\ -\frac{1}{x(1-x)} & \frac{1-x}{(1-y)^2} +
    \frac{x}{y^2} \end{array} \biggr ], \qquad \mbox{and} 
\qquad \nabla^2 g(x,y) =
\biggr [\begin{array}{cc} \frac{1-y}{(1-x)^2} + \frac{y}{x^2} &
    -\frac{1}{x(1-x)} \\ -\frac{1}{x(1-x)} & \frac{1}{y(1-y)}
    \end{array} \biggr ].
\end{eqnarray*}
Noting that $\derivone{f}{x} = -\derivone{g}{y};\; 
\derivtwo{g}{x}{y} = -\derivtwo{f}{x}{x};\; 
\derivtwo{f}{x}{y} = -\derivtwo{g}{y}{y}$, we
see that equation~\eqref{Eqn-Algebra-1} is equivalent to
\begin{eqnarray}
\label{Eqn-Algebra-Bounds}
2(\derivone{f}{x}\derivone{g}{x}\derivtwo{f}{y}{y} +
\derivone{f}{y}\derivone{g}{x}\derivtwo{g}{y}{y}) -
\derivone{f}{x}^2\derivtwo{g}{y}{y} +
\derivone{f}{y}^2\derivtwo{g}{x}{x} & \geq &
\frac{2}{g}(\derivone{f}{y}\derivone{g}{x} -
\derivone{f}{x}\derivone{g}{y})^2.
\end{eqnarray}
To simplify the algebra further, we shall make use of the inequality
$(\log t^2)^2 \leq (t-1/t)^2$, which is valid for any $t$. This
implies that
\[\derivone{f}{y}\derivone{g}{x} = (v-u)(1/v-1/u) \leq
\derivone{f}{x}\derivone{g}{y}= -(\log(u/v))^2 
= - \derivone{f}{x}^2 
= -\derivone{g}{y}^2\leq 0.\]
Thus, $-\derivone{f}{x}^2\derivtwo{g}{y}{y} \geq 
\derivone{f}{y}\derivone{g}{x}\derivtwo{g}{y}{y}$, and  
$\frac{2}{g}(\derivone{f}{y}\derivone{g}{x} - 
\derivone{f}{x}\derivone{g}{y})^2 \leq
\frac{2}{g}\derivone{f}{y}\derivone{g}{x}(
\derivone{f}{y}\derivone{g}{x} - \derivone{f}{x}\derivone{g}{y}).$
As a result, \eqref{Eqn-Algebra-Bounds} would follow if we can
show that
\[2(\derivone{f}{x}\derivone{g}{x}\derivtwo{f}{y}{y} +
\derivone{f}{y}\derivone{g}{x}\derivtwo{g}{y}{y}) +
\derivone{f}{y}\derivone{g}{x}\derivtwo{g}{y}{y} +
\derivone{f}{y}^2\derivtwo{g}{x}{x} \geq
\frac{2}{g}\derivone{f}{y}\derivone{g}{x}(
\derivone{f}{y}\derivone{g}{x} - \derivone{f}{x}\derivone{g}{y}).\]
For all $x \neq y$, we may divide both sides by $-\derivone{f}{y}(x,y)
\derivone{g}{x}(x,y) > 0$.  Consequently, it suffices to show that:
\begin{equation*}
-2\derivone{f}{x}\derivtwo{f}{y}{y}/\derivone{f}{y} -
\derivone{f}{y}\derivtwo{g}{x}{x}/\derivone{g}{x} -
3\derivtwo{g}{y}{y} \geq \frac{2}{g}(\derivone{f}{x}\derivone{g}{y} -
\derivone{g}{x}\derivone{f}{y}),
\end{equation*}
or, equivalently,
\begin{equation*}
2\log(u/v)\biggr (\frac{v}{u-1} + \frac{u}{1-v} \biggr ) + \biggr (
\frac{u}{1-x} + \frac{v}{x} \biggr ) - \frac{3}{y(1-y)} \geq
\frac{2}{g}\biggr (\frac{(u-v)^2}{uv} - (\log\frac{u}{v})^2 \biggr ),
\end{equation*}
or, equivalently,
\begin{equation}
\label{EqnFinalMess}
2\log(u/v)\frac{(u-v)(u+v-1)}{(u-1)(1-v)} +
\frac{(u-v)^2(u+v-4uv)}{uv(u-1)(1-v)}
\geq \frac{2}{g} \biggr (\frac{(u-v)^2}{uv} - (\log\frac{u}{v})^2 \biggr ).
\end{equation}
Due to the symmetry, it suffices to prove~\eqref{EqnFinalMess} for $x < y$.
In particular, we shall use the
following inequality for logarithm mean~\cite{Inequalities93},
which holds for $u\neq v$:
\[\frac{3}{2\sqrt{uv} + (u+v)/2} < \frac{\log u - \log v}{u-v}
< \frac{1}{(uv(u+v)/2)^{1/3}}.\]
We shall replace $\frac{\log(u/v)}{u-v}$ in~\eqref{EqnFinalMess} by
appropriate upper and lower bounds. In addition, we shall
also bound $g(x,y)$ from below, using the following argument.
When $x<y$, we have $u<1<v$, and
\begin{eqnarray*}
g(x,y) & = & y\log \frac{y}{x} + (1-y)\log \frac{1-y}{1-x}
> \frac{3y(y-x)}{2\sqrt{xy} + (x+y)/2} +
  \frac{(1-y)(x-y)}{[(1-x)(1-y)(1-(x+y)/2)]^{1/3}} \\
& = & \frac{3(1-v)(1-u)}{(u-v)(2\sqrt{u} + \frac{u+1}{2})}
+ \frac{(u-1)(1-v)}{(u-v)(v(v+1)/2)^{1/3}} > 0.
\end{eqnarray*}
\comment{
and if $x>y$,
\[g(x,y)  =  y \log \frac{1/x}{1/y} + (1-y) \log \frac{1-y}{1-x} 
> \frac {y(1/x-1/y)}{[1/(xy) (1/x+1/y)/2]^{1/3}} +
  \frac{3(1-y)(x-y)}{2\sqrt{(1-x)(1-y)} + (1-(x+y)/2)}\]
}
Let us denote this lower bound by $q(u,v)$.

Having got rid of the logarithm terms,~\eqref{EqnFinalMess} will 
hold if we can prove the following:
\[ \frac{6(u-v)^2(u+v-1)}{(2\sqrt{uv} + (u+v)/2)(u-1)(1-v)} + 
\frac{(u-v)^2(u+v-4uv)}{uv(u-1)(1-v)}
\geq\frac{2}{q(u,v)} \biggr (\frac{(u-v)^2}{uv} - \frac{9(u-v)^2}
{(2\sqrt{uv} + (u+v)/2)^2}\biggr ),\]
\comment{
or equivalently,
\begin{equation}
\frac{6(u+v-1)}{(2\sqrt{uv} + (u+v)/2)(u-1)(1-v)} + 
\frac{(u+v-4uv)}{uv(u-1)(1-v)}
\geq\frac{2}{q(u,v)} \biggr (\frac{1}{uv} - \frac{9}
{(2\sqrt{uv} + (u+v)/2)^2}\biggr ),
\end{equation}
}
or equivalently,
\begin{multline}
\label{VeryFinalMess}
\biggr ( \frac{6(u+v-1)}{(2\sqrt{uv} + (u+v)/2)} + 
\frac{(u+v-4uv)}{uv} \biggr )
\biggr ( \frac{3}{(v-u)(2\sqrt{u}+\frac{u+1}{2})} -
\frac{1}{(v-u)(v(v+1)/2)^{1/3}} \biggr ) \\
\geq 2 \biggr (\frac{1}{uv} - \frac{9}
{(2\sqrt{uv} + (u+v)/2)^2}\biggr ),
\end{multline}
which is equivalent to
\begin{multline}
\frac{(u+v-2\sqrt{uv})((u+v)/2 + 3\sqrt{uv} + 4uv)}
{(2\sqrt{uv} + (u+v)/2)uv}
\frac{3(v(v+1)/2)^{1/3} - (2\sqrt{u} + (u+1)/2)}
{(v-u)(2\sqrt{u} + (u+1)/2)(v(v+1)/2)^{1/3}} \\
\geq
\frac{(u+v - 2\sqrt{uv})((u+v)/2 + 5\sqrt{uv})}
{uv(2\sqrt{uv} + (u+v)/2)^2}
\end{multline}
and also equivalent to
\begin{multline}
\label{VeryVeryFinalMess}
((u+v)/2 + 2\sqrt{uv})((u+v)/2 + 3\sqrt{uv} + 4uv)
[3(v(v+1)/2)^{1/3} - (2\sqrt{u} + (u+1)/2)] \\ \geq
(2\sqrt{u} + (u+1)/2)(v(v+1)/2)^{1/3}((u+v)/2 + 5\sqrt{uv})(v-u).
\end{multline}

It can be checked by tedious but straightforward calculus that
inequality~\eqref{VeryVeryFinalMess} holds for any $u \leq 1 \leq v$, and equality
holds when $u=1=v$, i.e., $x=y$.

\mysection{Proof of Theorem~\ref{ThmPartialConv}}
\label{AppPartialConv}

Suppose that $\phi$ is not a likelihood ratio rule.
Then there exist positive $\Prob_0$-probability disjoint sets
$S_1,S_2,S_3$ such that for any $X_1 \in S_1, X_2 \in S_2, X_3 \in
S_3$,
\begin{subequations}
\label{Ineq-ratio}
\begin{eqnarray}
&&\phi(X_1) = \phi(X_3) = u_1 \\ &&\phi(X_2) = u_2 \neq u_1 \\
&&\frac{f^1(X_1)}{f^0(X_1)} < \frac{f^1(X_2)}{f^0(X_2)} <
\frac{f^1(X_3)}{f^0(X_3)}.
\end{eqnarray}
\end{subequations}
Define the probability of the quantiles as:
\begin{eqnarray*}
f^0(u_1) \defn \Prob_0(\phi(X) = u_1), & \mbox{and}& f^0(u_2) \defn
\Prob_0(\phi(X) = u_2), \\
f^1(u_1) \defn \Prob_1(\phi(X) = u_1), & \mbox{and} & f^1(u_2) \defn
\Prob_1(\phi(X) = u_2).
\end{eqnarray*}
Similarly, for the sets $S_1, S_2$ and $S_3$, we define
\begin{eqnarray*}
a_0 = f^0(S_1),\quad b_0 = f^0(S_2) & \mbox{and} & c_0 = f^0(S_3), \\
a_1 = f^1(S_1), \quad b_1 = f^1(S_2), & \mbox{and} & c_1 = f^1(S_3).
\end{eqnarray*}
Finally, let $p_0, p_1, q_0$ and $q_1$ denote the probability measures 
of the ``residuals'':
\begin{eqnarray*}
p_0 & = & f^0(u_2) - b_0, \qquad p_1 \; =\; f^1(u_2) - b_1, \\ q_0 & =
& f^0(u_1) - a_0 - c_0, \qquad q_1 \; =\; f^1(u_1) - a_1 - c_1.
\end{eqnarray*}
Note that we have $\frac{a_1}{a_0} < \frac{b_1}{b_0} <
\frac{c_1}{c_0}$.  In addition, the sets $S_1$ and $S_3$ were chosen
so that $\frac{a_1}{a_0} \leq \frac{q_1}{q_0} \leq \frac{c_1}{c_0}$.
From Lemma~\ref{Prop-PseudoLR}, there holds $\frac{p_1+b_1}{p_0+b_0} =
\frac{f^1(u_2)}{f^0(u_2)} \notin \biggr (\frac{a_1}{a_0},
\frac{c_1}{c_0} \biggr )$. We may assume without loss of generality
that $\frac{p_1+b_1}{p_0+b_0} \leq \frac{a_1}{a_0}$. Then, 
$\frac{p_1+b_1}{p_0 + b_0} < \frac{b_1}{b_0}$, so 
$\frac{p_1}{p_0} < \frac{p_1+b_1}{p_0 + b_0}$. Overall, we are
guaranteed to have the ordering
\begin{equation}
\label{EqnOrdering1}
\frac{p_1}{p_0} < \frac{p_1+b_1}{p_0+b_0} \leq
\frac{a_1}{a_0} < \frac{b_1}{b_0} < \frac{c_1}{c_0}.
\end{equation}

\newcommand{\setA}{\ensuremath{\mathcal{A}}}
\newcommand{\setB}{\ensuremath{\mathcal{B}}}
\newcommand{\setC}{\ensuremath{\mathcal{C}}}
\newcommand{\setP}{\ensuremath{\mathcal{P}}}
\newcommand{\setQ}{\ensuremath{\mathcal{Q}}}

Our strategy will be to modify the quantizer $\phi$ only for those
$X$ for which $\phi(X)$ takes the values $u_1$ or $u_2$, such that
the resulting quantizer is defined by a LLR-based threshold, and has 
a smaller (or equal) value of the corresponding cost $J_\phi^*$.  
For simplicity in notation, we use $\setA$ to denote the set with measures under
$\Prob_0$ and $\Prob_1$ equal to $a_0$ and $a_1$; the sets $\setB$,
$\setC$, $\setP$ and $\setQ$ are defined in an analogous manner.  We
begin by observing that we have either $\frac{a_1}{a_0} \leq
\frac{q_1+a_1}{q_0+a_0} < \frac{b_1}{b_0}$ or $\frac{b_1}{b_0} <
\frac{q_1+c_1}{q_0+c_0} \leq \frac{c_1}{c_0}$. Thus, in our subsequent
manipulation of sets, we always bundle $\setQ$ with either $\setA$ or
$\setC$ accordingly without changing the ordering of the probability
ratios.  Without loss of generality, then, we may disregard the
corresponding residual set corresponding to $\setQ$ in the analysis to
follow.

In the remainder of the proof, we shall show that either one of the
following two modifications of the quantizer $\phi$ will improve
(decrease) the sequential cost $J_\phi^*$: 
\begin{enumerate}
\item[(i)] Assign $\setA, \setB$ and $\setC$ to the same quantization
level $u_1$, and leave $\setP$ to the level $u_2$, or
\item[(ii)] Assign $\setP$, $\setA$ and $\setB$ to the same level
$u_2$, and leave $c$ to the level $u_1$.
\end{enumerate}
It is clear that this modified quantizer design respects the
likelihood ratio rule for the quantization indices $u_1$ and
$u_2$. By repeated application of this modification for every such
pair, we are guaranteed to arrive at a likelihood ratio quantizer that
is optimal, thereby completing the proof.

Let $a'_0,b'_0,c'_0,p'_0$ be normalized versions of $a_0, b_0, c_0,
p_0$, respectively (i.e., $a'_0 = a_0/(p_0+a_0+b_0+c_0)$, and so
on). Similarly, let $a'_1,b'_1,c'_1,p'_1$ be normalized versions of
$a_1, b_1, c_1, p_1$, respectively.  With this notation, we have the
relations

\begin{eqnarray*}
\asymKL{\phi} & = & \sum_{u\neq u_1,u_2}
f^0(u)\log\frac{f^0(u)}{f^1(u)} + (p_0 +
b_0)\log\frac{p_0+b_0}{p_1+b_1} + (a_0+c_0)\log \frac{a_0 +
c_0}{a_1+c_1} \\ & = & A_0 + (f^0(u_1)+f^0(u_2))\biggr((p'_0 +
b'_0)\log\frac{p'_0+b'_0} {p'_1+b'_1} + (a'_0+c'_0)\log \frac{a'_0 +
c'_0}{a'_1+c'_1}\biggr ) \\ & = & A_0 +
(f^0(u_1)+f^0(u_2))\extKL{p'+b'}{a'+c'}, \\
\asymLK{\phi} & = & \sum_{u\neq u_1,u_2} f^1(u)\log\frac{f^1(u)}{f^0(u)}
+ (p_1 + b_1)\log\frac{p_1+b_1}{p_0+b_0} + (a_1+c_1)\log
\frac{a_1 + c_1}{a_0+c_0} \\
%
%
& = & A_1 + (f^1(u_1)+f^1(u_2))\extLK{p'+b'}{a'+c'},
\end{eqnarray*}
where we define 
\begin{eqnarray*}
A_0 & \defn & \sum_{u\neq u_1,u_2}
f^0(u)\log\frac{f^0(u)}{f^1(u)} +
(f^0(u_1)+f^0(u_2))\log\frac{f^0(u_1)+f^0(u_2)}{f^1(u_1)+f^1(u_2)}
\geq 0, \\ 
A_1 & \defn & \sum_{u\neq u_1,u_2}
f^1(u)\log\frac{f^1(u)}{f^0(u)} +
(f^1(u_1)+f^1(u_2))\log\frac{f^1(u_1)+f^1(u_2)}{f^0(u_1)+f^0(u_2)}
\geq 0
\end{eqnarray*}
due to the non-negativity of the KL divergences.

Note that from~\eqref{EqnOrdering1} we have
\[\frac{p'_1}{p'_0} < \frac{p'_1+b'_1}{p'_0+b'_0} \leq
\frac{a'_1}{a'_0} < \frac{b'_1}{b'_0} < \frac{c'_1}{c'_0},\] in
addition to the normalization constraints that $p'_0 + a'_0 + b'_0 +
c'_0 = p'_1 + a'_1 + b'_1 + c'_1 = 1$.  It follows that
$\frac{p'_1+b'_1}{p'_0+b'_0} < \frac{p'_1 + a'_1 + b'_1 + c'_1}
{p'_0 + a'_0 + b'_0 + c'_0} = 1$. 

Let us consider varying the values of $a'_1,b'_1$, while fixing all
other variables and ensuring that all the above constraints
hold. Then, $a'_1 + b'_1$ is constant, and both $\extKL{p'+b'}{a'+c'}$
and $\extLK{p'+b'}{a'+c'}$ increase as $b'_1$ decreases and $a'_1$
increases. In other words, if we define $a''_0 = a'_0$, $b''_0 = b'_0$
and $a''_1$ and $b''_1$ such that
\[\frac{a''_1}{a'_0} = \frac{b''_1}{b'_0} = 
\frac{1-p'_1 - c'_1}{1-p'_0-c'_0},\]
then we have
\begin{equation}
\label{Eqn-Step-1}
\extKL{p'+b'}{a'+c'} \leq \extKL{p'+b''}{a''+c'} \;\mbox{and}\;
  \extLK{p'+b'}{a'+c'} \leq \extLK{p'+b''}{p''+c'}.
\end{equation}


Now note that vector $(b''_0, b''_1)$ in $\real^2$ is a convex
combination of $(0,0)$ and $(a''_0+b''_0, a''_1+b''_1)$.  It follows
that $(p'_0+b''_0,p'_1+b''_1)$ is a convex combination of
$(p'_0,p'_1)$ and $(p'_0+a''_0+b''_0, p'_1+a''_1+b''_1) = (p'_0 + a'_0
+ b'_0, p'_1 + a'_1 + b'_1)$. 

By ~\eqref{Eqn-Step-1}, we obtain:
\begin{eqnarray*}
G_\phi 
& = & \frac{\pi^0}{A_0 +
(f^0(u_1)+f^0(u_2))\extKL{p'+b'}{a'+c'}} + \frac{\pi^1}{A_1 +
(f^1(u_1)+f^1(u_2))\extLK{p'+b'}{a'+c'}} \\ & \geq & \frac{\pi^0}{A_0
+ (f^0(u_1)+f^0(u_2))\extKL{p'+b''}{a''+c'}} + \frac{\pi^1}{A_1 +
(f^1(u_1)+f^1(u_2))\extLK{p'+b''}{a''+c'}} \\ 
& = &
\frac{\pi^0}{A_0 + (f^0(u_1)+f^0(u_2))D(p'_0+b''_0, p'_1 + b''_1)} + 
\frac{\pi^1}{A_1 + (f^1(u_1)+f^1(u_2))D(p'_1+ b''_1, p'_0 + b''_0)}
\end{eqnarray*}
Applying the quasiconcavity result in Lemma~\ref{Lem-Quasiconcave}:
\begin{eqnarray*} 
G_\phi & \geq & 
\min \biggr \{
\frac{\pi^0}{A_0 + (f^0(u_1)+f^0(u_2))D(p'_0, p'_1)} +
\frac{\pi^1}{A_1 + (f^1(u_1)+f^1(u_2))D(p'_1, p'_0)}, \\ &&
\frac{\pi^0}{A_0 + (f^0(u_1)+f^0(u_2))D(p'_0+a'_0+b'_0, p'_1 + a'_1 + b'_1)} + \\ &&
\frac{\pi^1}{A_1 + (f^1(u_1)+f^1(u_2))D(p'_1+a'_1+b'_1, p'_0 + a'_0 + b'_0)} \biggr \}.
\end{eqnarray*}
But the two arguments of the minimum are the sequential cost
coefficient corresponding to the two possible modifications 
of $\phi$. Hence, the proof is complete.



\bibliography{sequential} \end{document}